\documentclass[a4paper, 10pt, DIV10, headinclude=false, footinclude=false]{scrartcl}

\usepackage[utf8]{inputenc}
\usepackage[T1]{fontenc}
\usepackage[english]{babel}

\usepackage[pdftex]{graphicx}
\usepackage{subcaption}
\usepackage{latexsym}
\usepackage{amsmath,amssymb,amsthm}
\usepackage{esint}
\usepackage{aligned-overset}
\usepackage{multirow}
\usepackage{bigdelim}
\usepackage{empheq} 
\usepackage{tikz}
\usepackage{booktabs}
\usetikzlibrary{positioning,shapes.multipart}

\newtheorem{Conclusion}{Conclusion}	
\theoremstyle{definition}

\theoremstyle{remark}	  
                  
\numberwithin{equation}{section} 
\numberwithin{figure}{section}

\usepackage{hyperref}
\hypersetup{
  pdffitwindow=false,
  pdfhighlight=/O,
  pdfnewwindow,
  colorlinks=true,
  citecolor=red,      
  linkcolor=blue,        
  menucolor=blue,       
  urlcolor=blue,          
  pdfpagemode=UseOutlines,
  bookmarksnumbered=true,
  linktocpage,
  pdfkeywords={},
  pdfcreator={pdflatex},
  pdfproducer={LaTeX with hyperref}
}

\newcommand{\R}{\mathbb{R}}
\newcommand{\N}{\mathbb{N}}

\newcommand{\bfE}{\mathbf{E}}
\newcommand{\bfH}{\mathbf{H}}
\newcommand{\bfD}{\mathbf{D}}
\newcommand{\bfB}{\mathbf{B}}
\newcommand{\bfu}{\mathbf{u}}
\newcommand{\bfv}{\mathbf{v}}
\newcommand{\bfw}{\mathbf{w}}

\newcommand{\el}{\mathrm{e}}
\newcommand{\ma}{\mathrm{m}}

\newcommand{\ve}{\bfv^{\el}}
\newcommand{\vm}{\bfv^{\ma}}

\newcommand{\obfu}{\overline{\bfu}}
\newcommand{\tbfu}{\widetilde{\bfu}}
\newcommand{\sobfu}[1]{\obfu^{(#1)}}
\newcommand{\sbfu}[1]{\bfu^{(#1)}}

\newcommand{\sobfD}[1]{\obfD^{(#1)}}
\newcommand{\sbfD}[1]{\bfD^{(#1)}}
\newcommand{\sobfB}[1]{\obfB^{(#1)}}
\newcommand{\sbfB}[1]{\bfB^{(#1)}}

\newcommand{\obfD}{\overline{\bfD}}
\newcommand{\tbfD}{\widetilde{\bfD}}
\newcommand{\obfB}{\overline{\bfB}}
\newcommand{\tbfB}{\widetilde{\bfB}}

\newcommand{\cor}{\mathrm{cor}}

\DeclareMathOperator{\curl}{\operatorname{curl}}
\DeclareMathOperator{\Div}{\operatorname{div}}

\newenvironment{abstr}[1]{ \vspace{.05in}\footnotesize
	\parindent .2in
	{\upshape\bfseries #1. }\ignorespaces}{\par\vspace{.1in}}
\newenvironment{Abstract}{\begin{abstr}{Abstract}}{\end{abstr}}
\newenvironment{keywords}{\begin{abstr}{Key words}}{\end{abstr}}
\newenvironment{AMS}{\begin{abstr}{AMS subject classifications}}{\end{abstr}}

\allowdisplaybreaks[4]

\begin{document}
	\title{Local and nonlocal homogenization of wave propagation in time-varying media%
		\thanks{Funded by the Deutsche Forschungsgemeinschaft (DFG, German Research Foundation) – Project-ID 258734477 – SFB 1173 and under Germany's Excellence Strategy – EXC-2047/1 – 390685813.}
	}
	\author{Christian D\"oding\footnotemark[2] \and Barbara Verf\"urth\footnotemark[2]}	
	
	\date{}
	\maketitle
	
	\renewcommand{\thefootnote}{\fnsymbol{footnote}}
	\footnotetext[2]{Institut für Numerische Simulation, Universit\"at Bonn, Friedrich-Hirzebruch-Allee 7, D-53115 Bonn, Germany}
	\renewcommand{\thefootnote}{\arabic{footnote}}
	
	\begin{Abstract} 
		Temporal metamaterials are artificially manufactured materials with time-dependent material properties that exhibit interesting phenomena when waves propagate through them. The propagation of electromagnetic waves in such time-varying dielectric media is governed by Maxwell's equations, which lead to wave equations with temporal highly oscillatory coefficients for the electric and magnetic fields. In this study, we analyze the effective behavior of electromagnetic fields in time-varying metamaterials using a formal two-scale asymptotic expansion. We provide a mathematical derivation of the effective equations for the leading-order homogenized solution, as well as for the first- and second-order corrections of the effective solution. While the effective solution and the first-order correction are governed by local material laws, we reveal a nonlocal constitutive relation for the second-order corrections. Special attention is also paid to temporal interface conditions through initial values of the homogenized equations. The results provide a mathematically justified framework for the effective description of wave-type equations of time-varying media, applicable to models in optics, elasticity, and acoustics.
	\end{Abstract}
	
	

\section{Introduction}

Metamaterials are artificially structured media with spatially and/or temporally varying material properties that enable wave phenomena beyond those found in natural materials, such as amplification, negative refraction indices, non-reciprocity, and cloaking. In recent years, temporal and spatiotemporal metamaterials have attracted significant attention due to their intriguing physical behavior and potential applications in optics, acoustics, and elasticity. We refer to the seminal works \cite{AFKR18, Cai10, Caloz20a, Caloz20, Cui10, Engheta21, Engheta06, Galiffi22, Mil23, SSB19} for a comprehensive introduction and overview of the theoretical foundations, design, and applications of spatiotemporal metamaterials from both physical and mathematical perspectives. The mathematical analysis of temporal metamaterials, in particular, has become an active area of research in recent decades. This field ranges from Morgenthaler’s early work on time-varying media \cite{Mor1958} to more recent studies in time-harmonic regimes \cite{Ammari25,Fante71,Garnier21,Gaxiola-Luna21,Torrent20,Xiao14,Zurita-Sanchez09}. The present work contributes to this research field from the standpoint of mathematical homogenization theory. As we elaborate below, it is closely related to studies on spatiotemporal metamaterials \cite{Dal25, FishChen01,Huidobro21,PaEn2020,RCG22,Touboul24a,Touboul24b,WautierGuzina15,WeZhWuLi2021}, in which different but conceptually related models of modulated media have been analyzed. Before discussing these connections, we first outline our model setting, describe the homogenization approach we exploit, and summarize our main contributions to the theory of homogenization for temporal metamaterials. \\
In this work we consider higher-order homogenization of temporally modulated electromagnetic metamaterials characterized by a spatially homogeneous but temporally heterogeneous permittivity. Electromagnetic wave propagation in such media is governed by Maxwell’s equations, which give rise to wave-type equations with time-dependent coefficients for the electric and magnetic fields. Using a formal two-scale asymptotic expansion, we derive effective wave equations for both fields. In addition to the leading-order homogenized equations, we carry out higher-order homogenization by deriving governing equations for the first- and second-order correctors of the effective solution. These higher-order terms yield refined approximations of the effective description and capture nonlocal physical effects as light propagates through the medium. We observe that the homogenized equations differ for the electric and magnetic fields due to the distinct placement of the time-dependent coefficients in their respective wave equations. Furthermore, a key finding of this study is the derivation of initial conditions that require careful treatment and introduce temporal interface conditions for the homogenized equations. Recasting our results within the Maxwell framework leads to both local and nonlocal effective material laws. Nonlocality means that the corresponding constitutive relation involves derivatives of the electric field, aligning with the literature, c.f. \cite{GMPR20,MKSPR18,RCG22}. While our primary focus is electromagnetics, the generality of our approach extends naturally to other wave systems in temporally modulated media, such as acoustics and elasticity.

Next, we discuss how our analysis generalizes and unifies results from previous studies, contributing to the existing literature. The work \cite{RCG22} addresses a case corresponding to the electric field setting considered here and derives second-order effective equations in the time-frequency domain, which are consistent with our results in physical space. Our contribution complements theirs by providing a general derivation of the homogenized equations, adding initial conditions, and extending the analysis to the magnetic field case. Higher-order spatial homogenization for wave-type equations has been studied in \cite{FishChen01,WautierGuzina15}. However, these settings do not cover the time-dependent case considered in our study, since the coefficients vary in space rather than in time. Although a formal interchange of space and time leads to similar homogenized equations, this approach does not cover initial values and temporal interface conditions or the case of space dimension larger than one. Spatiotemporal metamaterials with modulations along traveling-wave trajectories were studied at the effective level in \cite{Huidobro21}, and higher-order homogenization was subsequently investigated in \cite{Touboul24a,Touboul24b}. Roughly speaking, our setting belongs to case of infinite modulation velocity in \cite{Touboul24a, Touboul24b} which is excluded therein.

At leading order, an effective description for purely time-dependent dielectric composites was derived in \cite{PaEn2020}, where the authors obtain the harmonic mean as the effective permittivity for a two-component system with piecewise constant modulations. Their setting also corresponds to our electric field case. Although our framework assumes smooth modulations, it can be applied to the composite setting of \cite{PaEn2020}. In contrast, \cite{WeZhWuLi2021} examines time-varying acoustic media and finds that the arithmetic mean is the effective coefficient. At first, this seems to contradict \cite{PaEn2020}, but in \cite{WeZhWuLi2021}, the authors consider time-varying compressibility. Via the derivation of the wave equation, this corresponds to a different placement of the time-varying coefficient, which leads to a different homogenized coefficient, as highlighted in our analysis.
Finally, we note that spatiotemporal metamaterials were recently studied in \cite{Dal25} using a space-harmonic approach, which in certain cases allows for the derivation of exact solutions to the corresponding Maxwell system.

Summarizing, our findings provide a mathematically justified framework for simulating wave propagation in time-varying media. Our results lend theoretical support to prior observations and extend them. \\

\textbf{Outline.} 
The structure of the paper is as follows: In Section \ref{sec:setting}, we formulate the problem by presenting Maxwell’s equations and the associated wave-type equations for the electric and magnetic fields. Section \ref{sec:asympotic} introduces the formal two-scale asymptotic expansion that forms the basis of our analysis. First, we homogenize the electric field in Section \ref{subsec:hom:case2}, and then we present the analogous treatment for the magnetic field in Section \ref{subsec:hom:case1}. Section \ref{sec:Maxwell} connects our findings back to Maxwell’s equations, yielding the effective local and nonlocal material laws. Finally, Section \ref{sec:numerics} provides numerical simulations that demonstrate light propagation through time-varying media and validate our theoretical results.

\section{Problem derivation}\label{sec:setting}

Electromagnetic wave propagation in time-varying metamaterials is governed by Maxwell's equations in time domain $\R^3 \times \R_+$, reading in the absence of sources as
\begin{equation} \label{eq:Maxwell}
	\begin{aligned}
	\curl\bfE+\partial_t \bfB & = 0, \quad & \Div\bfD=0, \\
	\curl \bfH-\partial_t \bfD & =0, \quad & \Div\bfB=0.
	\end{aligned}
\end{equation}
Here $\bfD: \R^3 \times \R_+ \rightarrow \R^3$ is the electric displacement (induction), $\bfE: \R^3 \times \R_+ \rightarrow \R^3$ is the electric field, and $\bfB, \bfH: \R^3 \times \R_+ \rightarrow \R^3$ are the magnetic fields of the wave. This system is complemented with the linear constitutive relations
\begin{equation}
	\begin{split}
		\label{eq:material}
		\bfD(x,t)&=\varepsilon_\eta(t)\bfE(x,t),\qquad
		\bfB(x,t)=\mu(t)\bfH(x,t),
	\end{split}
\end{equation}
where the properties of the time-varying metamaterial are described by the permittivity $\varepsilon_\eta : \R_+ \rightarrow \R$, respectively the permeability $\mu: \R_+ \rightarrow \R$, that relates the electric components $\bfD$ and $\bfE$, respectively the magnetic components $\bfH$ and $\bfB$, of the electromagnetic wave. For simplicity, we only consider non-magnetic materials in the following and set $\mu \equiv 1$ so that $\bfB = \bfH$. In contrast, the electric permittivity is rapidly oscillating with respect to time on a microscopic length scale $0 < \eta \ll 1$, but spatially constant, modeling the time-varying property of the metamaterial. More specifically, we assume $\varepsilon_\eta$ to be periodic in time with period $\eta$. This corresponds to a time modulation of $\varepsilon_\eta$ with large modulation frequency $\omega_m = \eta^{-1}$. We denote by $\varepsilon$ the $1$-periodic blueprint so that
\begin{align*}
	\varepsilon_\eta(t) = \varepsilon(t/\eta), \quad 0 < \eta \ll 1.
\end{align*}
We emphasize that assuming $\varepsilon_\eta$ to be scalar and real corresponds to an idealized model of an isotropic, lossless metamaterial, since a non-zero imaginary part would account for material absorption. While allowing $\varepsilon_\eta$ to be complex-valued does not alter our analytical approach, and the resulting homogenized equations remain of the same form, the anisotropic case is much more involved and not covered by our approach. \\
Using the constitutive relations, Maxwell's equations reduce to a system for either $\bfE$ and $\bfH$ or $\bfD$ and $\bfB$. It is well known that one can combine the two time-dependent equations of \eqref{eq:Maxwell} to one single evolution for any of the four fields.
Incorporating the divergence constraints and using the identity $\curl \curl =-\Delta+\nabla\Div$ with the vector Laplacian $\Delta$, we deduce a wave-type equation for the selected component of the electromagnetic wave. We obtain three independent equations of different type for each of the fields $\bfB$, $\bfD$ and $\bfE$:
\begin{align}
		\label{eq:wave1}
		\partial_t(\varepsilon_\eta \partial_t\bfB)-\Delta \bfB= 0,\\
		\label{eq:wave2}
		\varepsilon_\eta \partial_{tt}\bfD - \Delta\bfD = 0,\\
		\label{eq:wave3}
		\partial_{tt}(\varepsilon_\eta \bfE )-\Delta\bfE =  0.
\end{align}
Obviously, the equations for the electric components \eqref{eq:wave2} and \eqref{eq:wave3} are strongly related and can be transferred into each other by the constitutive relations \eqref{eq:material} and it suffices to consider only \eqref{eq:wave2}. Thus, we end up with two different wave equations describing the time evolution of the magnetic field \eqref{eq:wave1} and of the electric displacement \eqref{eq:wave2}. \\
The system of Maxwell's equations is closed with the initial conditions
\begin{align*}
	\bfB(x,0) = \vm_0 \quad \text{and} \quad \bfD(x,0) = \ve_0
\end{align*}
for smooth initial vector fields $\vm_0: \R^3 \rightarrow \R^3$ and $\ve_0: \R^3 \rightarrow \R^3$. This determines the first initial conditions for the wave-type equations \eqref{eq:wave1} and \eqref{eq:wave2}. Both wave-type problems are closed with second initial conditions that can be derived from the temporal interface conditions of Morgenthaler \cite{Mor1958}. They require the curls of $\bfB$ and $\bfD$ to be continuous at $t=0$ leading to
\begin{equation} \label{eq:interface1}
	\varepsilon_\eta(0) \partial_t \bfB(x,0) = - \curl \bfD(x,0) = - \curl \ve_0
\end{equation}
and
\begin{equation} \label{eq:interface2}
	\partial_t \bfD(x,0) = \curl \bfB(x,0) = \curl \vm_0.
\end{equation}
Note that we discuss spatially unbounded media so that no spatial boundary conditions are required. Summarizing we obtain the following settings  -- PDE and initial conditions -- for the magnetic field and the electric displacement, where we now generalize to arbitrary space dimensions $d = 1, 2, 3$, arbitrary sufficiently smooth initial data $\ve_0$, $\ve_1$, $\vm_0$, $\vm_1: \R^d \rightarrow \R^d$ that may (or may not) satisfy the temporal interface conditions of Morgenthaler $\vm_1 = -\curl \ve_0$ and $\ve_1 = \curl \vm_0$.
\begin{description}
\item[Electric case:] The electric displacement $\bfD$ is a solution $\bfu_\eta$ to the initial value problem
\begin{equation} \label{eq:case2}
\begin{aligned}
	\varepsilon_\eta \partial_{tt} \bfu_\eta - \Delta \bfu_\eta & = 0 \quad && \text{in } \R^d \times \R_+, \\
	\bfu_\eta(\cdot,0) & = \ve_0 \quad && \text{in } \R^d, \\
	\partial_t \bfu_\eta(\cdot,0) & = \ve_1 \quad && \text{in } \R^d.
\end{aligned}
\end{equation}

\item[Magnetic case:] The magnetic field $\bfB$ is a solution $\bfu_\eta$ to the initial value problem
\begin{equation} \label{eq:case1}
\begin{aligned}
	\partial_t (\varepsilon_\eta \partial_t \bfu_\eta) - \Delta \bfu_\eta & = 0 \quad && \text{in } \R^d \times \R_+, \\
	\bfu_\eta(\cdot,0) & = \vm_0 \quad && \text{in } \R^d, \\
	\varepsilon_\eta(0) \partial_t \bfu_\eta(\cdot,0) & = \vm_1 \quad && \text{in } \R^d.
\end{aligned}
\end{equation}
\end{description}
Our goal is to derive an effective description of the solutions to \eqref{eq:case2} and \eqref{eq:case1}, that is, we consider the limit $\eta\to0$ of very fast time modulations of the material. Due to the time dependence of the permittivity $\varepsilon_\eta$ the two wave equations have different mathematical nature so that we expect also a different effective behavior and need to treat the two wave equations \eqref{eq:case2} and \eqref{eq:case1} separately.  We also derive formulas for the first and second order corrections of the effective behavior, which vanish as $\eta \rightarrow 0$, but lead to higher order approximations of the solutions in the case of large but finite modulation frequencies $\omega_m = \eta^{-1}$. In the latter case we strictly assume that the time modulations are much larger than the spatial modulations, i.e., $\omega_0 \ll \omega_m$, where $\omega_0$ denotes the spatial modulation frequency determined by the initial values. Therefore, in the asymptotic regime $\eta \rightarrow 0$ this condition is always satisfied, but one has to take care of it in particular simulations when the time modulation is finite.

\section{Homogenization via asymptotic expansion}\label{sec:asympotic}

In this section we introduce the formal two-scale asymptotic expansion of the solution $\bfu_\eta$ of \eqref{eq:case2} and \eqref{eq:case1}, respectively. We define the microscopic time variable $\tau = t/\eta$ and expand $\bfu_\eta$ into
\begin{equation}\label{eq:asymptoticexpansion}
	\mathbf u_\eta(x,t)=\mathbf u_0(x,t,\tau)+\eta\mathbf u_1(x,t,\tau)+\eta^2\mathbf u_2(x,t,\tau)+ \cdots = \sum_{j = 0}^k \eta^j \bfu_j(x,t,\tau) + \mathcal{O}(\eta^{k+1})
\end{equation}
for some $k \in \N$ and where each $\bfu_j$ is $1$-periodic w.r.t. $\tau$. With the introduced notation, the time derivative $\partial_t$ yields $\partial_t+\eta^{-1}\partial_\tau$ for any $\mathbf u_j$ due to the chain rule.
We plug-in this ansatz into the corresponding wave equations \eqref{eq:case2}, \eqref{eq:case1} and derive problems for the expansions functions $\mathbf u_0$, $\mathbf u_1$, $\mathbf u_2$. In particular, we derive problems for the macroscopic and microscopic components of the expansions functions $\bfu_j$ by the splitting
\begin{align} \label{eq:decomp}
	\bfu_j(x,t,\tau) = \obfu_j(x,t) + \tbfu_j(x,t,\tau), \quad j \in \N_0.
\end{align}
Here, $\obfu_j(x,t)$ is the macroscopic, $\tau$-independent component of the $k$-th expansion function (sometimes also called time-averaged part) and $\tbfu_j$ is the microscopic component which is $1$-periodic w.r.t. $\tau$ and has zeromean, i.e.,
\begin{align} \label{zero-mean}
	\int_0^1 \tbfu_j(x,t,\tau) d\tau = 0
\end{align}
for all $(x,t) \in \R^d \times \R_+$. The decomposition \eqref{eq:decomp} into macroscopic and microscopic components is motivated by separating the dominant, energy-carrying wave dynamics from the rapid, fine-scale oscillations induced by temporal modulation. Due to the zero-mean condition \eqref{zero-mean} these fine-scale contributions average out over each period. In the subsequent sections we first consider the equations for the electric displacement \eqref{eq:case2} (Section \ref{subsec:hom:case2}), which we frequently refer to as the electric case. Then we tackle the system \eqref{eq:case1} describing the propagation of the magnetic field (Section \ref{subsec:hom:case1}) which is referred to the magnetic case. \\
Clearly, for \eqref{eq:asymptoticexpansion} to hold we need to assume that the solution $\bfu_\eta$ of \eqref{eq:case1} and \eqref{eq:case2}, respectively, depends smoothly on $\eta$. This can be expected as long as the permittivity coefficient $\varepsilon_\eta$ and the initial data $\ve_0$, $\ve_1$ and $\vm_0$, $\vm_1$ are sufficiently smooth. Beyond assuming the formal expansion \eqref{eq:asymptoticexpansion}, let us briefly discuss the regularity requirements on the permittivity and initial conditions. It suffices to assume that the permittivity is bounded, i.e., $\varepsilon_\eta \in L^\infty(0,T;\R_+)$. Under this assumption, the homogenized coefficients are well-defined, and the subsequent formal derivations are justified. In the magnetic case, weak formulations in time are required to rigorously handle the equations. For the initial values $\ve_0, \ve_1$ and $\vm_0, \vm_1$ we assume regularity in Sobolev spaces $H^k(\Omega)$ with $k$ large enough to guarantee well-posedness of the subsequent homogenized wave equations in $L^2(0,T,H^1_0(
\Omega))$, c.f. \cite{Evans} \\
However, in this work we do not aim at a rigorous justification of the formal asymptotic expansion \eqref{eq:asymptoticexpansion}, but at the derivation and mathematical underpinning of the effective local and nonlocal material laws as electromagnetic waves propagate through time-varying metamaterials. In this sense, our results are purely formal, but several classical homogenization techniques can be used to rigorously justify our findings, e.g., the method of two-scale convergence \cite{All1992} or other related homogenization methods \cite{BeLiPa1978,CDG2019,CD99,Jikov94,Tartar09}. The techniques used therein have mainly been applied to spatial multiscale problems, but a suitable exchange of space and time might allow to treat time-dependent problems in a similar way. In fact, the media oscillate rapidly with respect to a one-dimensional variable, namely time, and homogenization in one dimension is an extensively studied area. Finally, for an application of these techniques, it is crucial that our media are constant in space and that the initial values do not contain any spatial multiscale features. This is expressed by the small spatial modulation frequencies imposed by the initial values in \eqref{eq:case2} and \eqref{eq:case1}, so that no coupling with the fast temporal scale can occur. 

\section{The electric case} \label{subsec:hom:case2}

In this section we consider the initial value problem \eqref{eq:case2} for the electric displacement and derive equations for the effective solution $\bfu_0$ and the higher order corrections $\bfu_j$, $j \ge 1$. The general strategy for the derivation relies on plugging the formal asymptotic expansion from \eqref{eq:asymptoticexpansion} into the PDE in \eqref{eq:case2} and then comparing powers in $\eta$. This leads to governing PDEs for the effective solution and the higher-order corrections. Similarly, the initial conditions from \eqref{eq:case2} are treated so that initial value problems for the effective solution and higher-order corrections are observed.

\subsection{Effective electric behavior} \label{sec:case2-eff}

We start with the effective behavior $\bfu_0$ in \eqref{eq:asymptoticexpansion}. First, we consider the PDE in \eqref{eq:case2} and discuss the initial conditions for $\bfu_0$ afterwards. Plugging in the formal asymptotic expansion \eqref{eq:asymptoticexpansion} into the PDE in \eqref{eq:case2} yields
\begin{align} \label{eq:case2-expansion}
\begin{split}
	0 = \varepsilon_\eta\partial_{tt}\bfu_\eta - \Delta \bfu_\eta &=\frac{1}{\eta^2}\varepsilon\partial_{\tau\tau}\bfu_0+\frac{1}{\eta}\varepsilon\Bigl(2\partial_{t\tau}\bfu_0+\partial_{\tau\tau}\bfu_1\Bigr) \\
	& \quad + \sum_{j = 0}^{k-2} \eta^j  \varepsilon \Bigl(\partial_{tt} \bfu_j+ 2 \partial_{t\tau} \bfu_{j+1}+ \partial_{\tau
\tau} \bfu_{j+2} - \varepsilon^{-1} \Delta \bfu_j \Bigr) +\mathcal O(\eta^{k-1}).
\end{split}
\end{align}
For the $\mathcal{O}(\eta^{-2})$-term we obtain
\begin{align*}
	\varepsilon\partial_{\tau\tau}\bfu_0 = 0.
\end{align*}
Since $\varepsilon$ is strictly positive, this implies that $\bfu_0$ is independent of $\tau$ and therefore its zero-mean component vanishes, $\tbfu_0 = 0$. Furthermore, the $\mathcal{O}(\eta^{-1})$-term yields
\begin{align} \label{eq:case2:eta1}
 	\varepsilon\partial_{\tau\tau}\bfu_1=0
\end{align}
and hence also $\tbfu_1 = 0$. The equation for the $\mathcal{O}(1)$-term reads after dividing by $\varepsilon$ as
\begin{equation} \label{eq:case2:eta0}
	\partial_{tt}\bfu_0- \varepsilon^{-1} \Delta \bfu_0 + \partial_{\tau\tau}\bfu_2 = 0.
\end{equation}
Integration over $(0,1)$ w.r.t. $\tau$ and using the periodicity of $\bfu_2$ yields
\begin{align*}
	\int_0^1 \partial_{tt}\bfu_0 - \varepsilon^{-1} \Delta \bfu_0 \, d\tau=0.
\end{align*}
Since $\bfu_0$ is independent of $\tau$ we conclude that the effective behavior $\bfu_0$ is given as a solution to the homogenized equation
\begin{equation}\label{eq:case2:homeq}
	\varepsilon_{\hom}\partial_{tt}\bfu_0-\Delta\bfu_0=0
\end{equation}
with the homogenized coefficient
\begin{equation}\label{eq:case1:harmav}
	\varepsilon_{\hom}=\Bigl(\int_0^1 \varepsilon^{-1} d\tau\Bigr)^{-1}.
\end{equation}
In other words, the effective permittivity in this case is the harmonic average of the time-modulated permittivity over one time period. This is exactly the well-known formula in one-dimensional spatial homogenization, see e.g. \cite{BeLiPa1978, CD99, Jikov94}. \\
Next we discuss the initial conditions for the effective equation \eqref{eq:case2:homeq}. The first initial condition, $\bfu_\eta(\cdot,0) = \ve_0$ in \eqref{eq:case2}, directly implies $\bfu_0(\cdot,0) = \ve_0$ (and $\bfu_1(x,0) = \bfu_2(x,0,0) = 0$). The second initial condition $\partial_t \bfu_\eta(\cdot,0) =\ve_1$ in \eqref{eq:case2} yields $\partial_t \bfu_0(\cdot,0) = \ve_1$, since $\bfu_1$ is independent of $\tau$. Our observations for the effective behavior $\bfu_0$ are summarized in the following conclusion.

\begin{Conclusion} \label{con:electric_hom}
Suppose the asymptotic expansion \eqref{eq:asymptoticexpansion} holds for the solution $\bfu_\eta$ of the initial value problem \eqref{eq:case2} up to order $k = 2$. Then $\tbfu_0 = 0$ in $\R^d \times \R_+$ and the effective solution $\bfu_0 = \obfu_0$ solves the homogenized initial value problem
\begin{equation} \label{eq:case2-homogenized}
\begin{aligned}
	\varepsilon_{\hom} \partial_{tt} \bfu_0 - \Delta \bfu_0 & = 0 \quad && \text{in } \R^d \times \R_+, \\
	\bfu_0(\cdot,0) & = \ve_0 \quad && \text{in } \R^d, \\
	\partial_t \bfu_0(\cdot,0) & = \ve_1 \quad && \text{in } \R^d
\end{aligned}
\end{equation}
with the homogenized coefficient $\varepsilon_{\hom }$ from \eqref{eq:case1:harmav}.
\end{Conclusion}

\subsection{Electric first order correction}

We continue with the problem derivation for the first order corrector $\bfu_1$ for which we have already seen that it is independent of $\tau$; cf. \eqref{eq:case2:eta1}. Recalling the expansion from \eqref{eq:case2-expansion}, we see that the $\mathcal{O}(\eta)$-term is
\begin{align} \label{electricetaterm}
\partial_{\tau\tau}\bfu_3+2\partial_{t\tau}\bfu_2+\partial_{tt}\bfu_1-\varepsilon^{-1}\Delta \bfu_1=0.
\end{align}
Again, integrating over $(0,1)$ w.r.t. $\tau$ and taking the periodicity of $\bfu_2$ and $\bfu_3$ into account yields, similar to the derivation of the homogenized equation \eqref{eq:case2:homeq}, that
\begin{equation*}
	\varepsilon_{\hom} \partial_{tt}\bfu_1 - \Delta \bfu_1=0.
\end{equation*}
Hence, $\bfu_1$ also solves the homogenized equation as $\bfu_0$ does, but the initial conditions differ as we now discuss. The first initial condition in \eqref{eq:case2} directly implies $\bfu_1(x,0) = 0$ and we again focus on the second initial condition $\partial_t \bfu_\eta(\cdot,0) = \bfv_1$ from \eqref{eq:case2}. Plugging in the asymptotic expansion \eqref{eq:asymptoticexpansion} into the initial condition we obtain for the $\mathcal{O}(\eta)$-term
\begin{align*}
	\partial_t \bfu_1(x,0) + \partial_\tau \bfu_2(x,0,0) = 0.
\end{align*}
We can calculate $\partial_\tau \bfu_2(x,0,0)$ explicitly. For this, we recall the $\mathcal{O}(\eta)$-term of the asymptotic expansion in the PDE, see \eqref{eq:case2:eta0}, and integrate it with respect to $\tau$ to obtain
\begin{align*}
	\partial_\tau \bfu_2 (\cdot,\cdot,\tau) & - \partial_\tau \bfu_2 (\cdot,\cdot,0) = \int_0^\tau \partial_{\tau\tau} \bfu_2 ds \\
	& = - \int_0^\tau \partial_{tt}\bfu_0-\varepsilon^{-1}\Delta \bfu_0\, ds = \Bigl(\int_0^\tau\varepsilon^{-1}ds\Bigr)\Delta \bfu_0 -  \tau \partial_{tt}\bfu_0.
\end{align*}
Another integration over $(0,1)$ w.r.t. $\tau$ yields, using $\int_0^1 \partial_\tau \bfu_2 d\tau = 0$ due to periodicity,
\begin{align*}
\partial_t \bfu_1(x,0)=-\partial_\tau \bfu_2(x,0,0)& =\Bigl(\int_0^1\int_0^\tau \varepsilon^{-1}ds\, d\tau\Bigr)\Delta\bfu_0(x,0) - \frac{1}{2} \partial_{tt}\bfu_0(x,0) \\
& = \Bigl(\int_0^1\int_0^\tau \varepsilon^{-1}ds\, d\tau - \frac{1}{2 \varepsilon_{\hom} } \Bigr)\Delta \ve_0 \\
& = \Bigl(\int_0^1 (1-s) \, \varepsilon^{-1}(s) ds - \frac{1}{2 \varepsilon_{\hom} } \Bigr)\Delta \ve_0 \\
& = - \chi_0 \varepsilon_{\hom}^{-1} \Delta \ve_0,
\end{align*}
where we used that $\bfu_0$ is a smooth at $t = 0$ and set
\begin{align} \label{chi0}
	\chi_0:= \frac{1}{2} - \varepsilon_{\hom} \int_0^1 \int_0^\tau \varepsilon^{-1} ds d\tau.
\end{align}
We summarize our findings in the following conclusion.
\begin{Conclusion} \label{con:electric_cor1}
Suppose the asymptotic expansion \eqref{eq:asymptoticexpansion} holds for the solution $\bfu_\eta$ of the initial value problem \eqref{eq:case2} up to order $k = 3$. Then $\tbfu_1 = 0$ in $\R^d \times \R_+$ and the first order correction $\bfu_1 = \obfu_1$ solves the homogenized initial value problem
\begin{equation} \label{eq:case2-firstcorrectionsystem}
\begin{aligned}
	\varepsilon_{\hom} \partial_{tt} \bfu_1 - \Delta \bfu_1 & = 0 \quad && \text{in } \R^d \times \R_+, \\
	\bfu_1(\cdot,0) & = 0 \quad && \text{in } \R^d, \\
	\partial_t \bfu_1(\cdot,0) & = - \chi_0 \varepsilon_{\hom}^{-1} \Delta \ve_0 \quad && \text{in } \R^d
\end{aligned}
\end{equation}
with the homogenized coefficient $\varepsilon_{\hom }$ from \eqref{eq:case1:harmav}.
\end{Conclusion}
We emphasize that if the permittivity blueprint $\varepsilon$ satisfies
\begin{equation}\label{eq:condeps}
	\int_0^1 \left( s - \frac{1}{2} \right) \varepsilon^{-1}(s) ds = 0
\end{equation}
it follows $\chi_0 = 0$. In this case, we have $\bfu_1 \equiv 0$ on $\R^d \times \R_+$ and the effective solution $\bfu_0$ is expected to be a second order approximation of $\bfu_\eta$ in $\eta$. However, this is clearly a specific condition on the time variation of $\varepsilon$ and belongs to a degenerate case.

\subsection{Electric higher order corrections}
We next investigate how higher order terms in the asymptotic expansion can lead to refined effective material laws, especially for finite $\eta$ and derive a system of equations for the second order corrector $\bfu_2$ for the electric field in \eqref{eq:case2}. As we will see the resulting effective equations involve higher order derivatives. Recall, such higher-order effects are called nonlocal in line with the physics literature. \\
From \eqref{eq:case2:eta0}, we see that
\begin{equation*}
	\partial_{\tau\tau}\bfu_2 = \varepsilon^{-1}\Delta \bfu_0-\partial_{tt}\bfu_0 = (\varepsilon_{\hom} \varepsilon^{-1} - 1)\partial_{tt} \bfu_0,
\end{equation*}
where we used that $\bfu_0$ solves the homogenized equation. Thus we can write $\bfu_2$ as
\begin{align*}
	\bfu_2(x,t,\tau) = \theta(\tau) \partial_{tt} \bfu_0 + \overline{\bfu}_2(x,t)
\end{align*}
so that $\tbfu_2(x,t,\tau) = \theta(\tau) \partial_{tt}\bfu_0(x,t)$ and where $\theta$ is the unique $1$-periodic solution of the cell-problem
\begin{align} \label{electric_cell_problem}
	\partial_{\tau\tau} \theta = \varepsilon_{\hom }\varepsilon^{-1} - 1, \quad \int_0^1 \theta \, d\tau = 0.
\end{align}
We can calculate $\theta$ explicitly as
\begin{align*}
	\theta(\tau) = \theta_0 - \frac{1}{2} (\tau^2 - \tau) + \varepsilon_{\hom} \Big(\int_0^\tau \int_0^s \varepsilon^{-1} drds - \tau \int_0^1 \int_0^s \varepsilon^{-1} drds \Big), 
\end{align*}
where
\begin{equation} \label{eq:theta_0}
	\theta_0 = -\frac{1}{12} - \varepsilon_{\hom} \Big(\int_0^1 \int_0^\tau \int_0^s \varepsilon^{-1} drdsd\tau - \frac{1}{2} \int_0^1 \int_0^s \varepsilon^{-1} drds \Big).
\end{equation}
Next we derive the equation for $\overline{\bfu}_2$. The $\mathcal{O}(\eta^2)$-term in the asymptotic expansion gives
\begin{equation} \label{eq:comp1}
	\partial_{\tau\tau}\bfu_4+2\partial_{t\tau}\bfu_3+\partial_{tt}\bfu_2-\varepsilon^{-1}\Delta \bfu_2=0.
\end{equation}
Integrating \eqref{eq:comp1} over $(0,1)$ w.r.t. $\tau$  and using the periodicity of $\bfu_3$ and $\bfu_4$ yields
\begin{align} \label{eq:comp2}
	0 = \int_0^1\partial_{tt}\bfu_2 - \varepsilon^{-1}\Delta \bfu_2 \, d\tau = \int_0^1 \partial_{tt} \overline{\bfu}_2 - \varepsilon^{-1} \Delta \overline{\bfu}_2 - \varepsilon^{-1} \theta \Delta \partial_{tt} \bfu_0 \, d\tau .
\end{align}
Replacing $\Delta \partial_{tt} \bfu_0$ by $\varepsilon_{\hom}^{-1} \Delta^2 \bfu_0$ in \eqref{eq:comp2} due to \eqref{eq:case2-homogenized}, this yields the homogenized equation with a (nonlocal) correction
\begin{align*}
	\varepsilon_{\hom} \partial_{tt} \overline{\bfu}_2 - \Delta \overline{\bfu}_2 = \varepsilon_{\cor} \Delta^2 \bfu_0,
\end{align*}
where
\begin{align} \label{electric_correction_coefficient}
	\varepsilon_{\cor} = \int_0^1 \varepsilon^{-1} \theta \, d\tau
\end{align}
is the correction coefficient of the material. Next we discuss the initial conditions. The first initial condition $\bfu_\eta(\cdot,0) = \ve_0$ yields $\bfu_2(x,0,0) = 0$ which implies $\overline{\bfu}_2(\cdot,0) = - \theta_0 \partial_{tt} \bfu_0(\cdot,0) = - \varepsilon_{\hom}^{-1} \theta_0 \Delta \ve_0$ where we again used that $\bfu_0$ is a smooth solution of \eqref{eq:case2-homogenized}. Plugging in the asymptotic expansion into the second initial condition $\partial_t \bfu_\eta(\cdot,0) = \ve_1$ yields for the $\mathcal{O}(\eta)$-term
\begin{align*}
	0 = \partial_t \bfu_2(x,0,0) + \partial_\tau \bfu_3(x,0,0).
\end{align*}
Further, recalling the $\mathcal{O}(\eta)$-term from \eqref{electricetaterm} we obtain after integration w.r.t. $\tau$
\begin{align*}
	\partial_\tau \bfu_3(x,0,\tau) - \partial_\tau \bfu_3(x,0,0) & = \int_0^\tau \partial_{\tau \tau} \bfu_3(x,0,s) \, ds \\
	& = - \int_0^\tau 2\partial_{t \tau} \bfu_2 (x,0,s) + \partial_{tt} \bfu_1(x,0) - \varepsilon^{-1}(s) \Delta \bfu_1(x,0) \, ds \\
	& = - \int_0^\tau 2\partial_\tau \theta(s) \partial_{ttt} \bfu_0(x,0) + \partial_{tt} \bfu_1(x,0) - \varepsilon^{-1}(s) \Delta \bfu_1(x,0) \, ds \\
	& = - \int_0^\tau 2\partial_\tau \theta(s) \partial_{ttt} \bfu_0(x,0) \, ds,
\end{align*}
where we used that $\Delta \bfu_1(x,0) = 0$ and $\partial_{tt} \bfu_1(x,0) = \varepsilon_{\hom}^{-1} \Delta \bfu_1(x,0) = 0$. Another integration over $(0,1)$ w.r.t. $\tau$ yields
\begin{align*}
	-\partial_\tau \bfu_3(x,0,0) = -  2 \int_0^1 \int_0^\tau \partial_\tau \theta \, dsd\tau \, \partial_{ttt} \bfu_0(x,0) = 2\theta_0 \partial_{ttt} \bfu_0(x,0).
\end{align*}
This yields
\begin{align*}
	\partial_t \overline{\bfu}_2(x,0) = \partial_t \bfu_2(x,0,0) - \partial_t \tbfu_2(x,0,0) = \theta_0 \partial_{ttt} \bfu_0(x,0) = \varepsilon_{\hom}^{-1} \theta_0 \Delta \ve_1(x).
\end{align*}

\begin{Conclusion} \label{con:electric_cor2}
Suppose the asymptotic expansion \eqref{eq:asymptoticexpansion} holds for the solution $\bfu_\eta$ of the initial value problem \eqref{eq:case2} up to order $k = 4$. Then the second order correction $\bfu_2 = \obfu_2 + \tbfu_2$ is given by
\begin{equation} \label{eq:case1-correction2_tilde}
	\tbfu_2(x,t,\tau) = \theta(\tau) \partial_{tt} \bfu_0(x,t)
\end{equation}
with the solutions $\theta$ of the cell problem \eqref{electric_cell_problem}. Furthermore, with the homogenized coefficient $\varepsilon_{\hom}$ from \eqref{eq:case1:harmav} and the correction coefficient $\varepsilon_{\cor}$ from \eqref{electric_correction_coefficient}, $\overline{\bfu}_2$ is given by the solution to the second order corrector problem
\begin{equation} \label{eq:case2-correction2}
\begin{aligned}
	\varepsilon_{\hom} \partial_{tt} \overline{\bfu}_2 - \Delta \overline{\bfu}_2 & = \varepsilon_{\cor} \Delta^2 \bfu_0 \quad && \text{in } \R^d \times \R_+, \\
	\overline{\bfu}_2(\cdot,0) & = - \varepsilon_{\hom}^{-1} \theta_0 \Delta \ve_0 \quad && \text{in } \R^d, \\
	\partial_t \overline{\bfu}_2(\cdot,0) & = \varepsilon_{\hom}^{-1} \theta_0 \Delta \ve_1 \quad && \text{in } \R^d.
\end{aligned}
\end{equation}
\end{Conclusion}
Let us note here that computing $\obfu_2$ and clearly also $\bfu_2$ requires to first solve for $\bfu_0$ via \eqref{eq:case2-homogenized}. We emphasize that $\varepsilon_{\cor}$ is negative, since \eqref{electric_cell_problem}, the periodicity of $\theta$, and integration by parts implies
\begin{align*}
	\varepsilon_\cor = \varepsilon_{\hom}^{-1} \int_0^1 \theta (\partial_{\tau \tau} \theta + 1) d\tau = \varepsilon_{\hom}^{-1} \int_0^1 \theta \partial_{\tau \tau} \theta d\tau = - \varepsilon_{\hom}^{-1} \int_0^1 (\partial_\tau \theta)^2 d\tau < 0.
\end{align*}
Therefore, equation \eqref{eq:case2-correction2}, sometimes called a ``good'' Boussinesq equation, is well-posed (cf. \cite{Bona88,Linares93}). Furthermore, it is worth to mention that the replacement of $\Delta \partial_{tt} \bfu_0$ by $\varepsilon_{\hom}^{-1} \Delta^2 \bfu_0$ is somehow arbitrary but favorable in our case. A more general replacement is $\Delta \partial_{tt} \bfu_0 =  q \varepsilon_{\hom}^{-1} \Delta^2 \bfu_0 + (1-q) \Delta \partial_{tt} \bfu_0$, providing an additional degree of freedom $q \in (0,1)$. However, to obtain a classical Boussinesq-type equation involving only higher order spatial derivatives we set $q = 1$.

\subsection{Pure macroscopic description of the electric field}
\label{sec:pure_marco_electric}

In contrast to solving for the full components $\bfu_j$, which in particular includes potential microscopic information via $\tbfu_j$, one is usually only interested in the purely macroscopic effective behavior of the propagating wave. We can express this macroscopic behavior by the quantities
\begin{align} \label{macro_quantities}
	\sobfu{k}(x,t) = \sum_{j = 0}^k \eta^j \obfu_j(x,t), \quad k \in \N. 
\end{align}
Due to our previous derivations we have $\sobfu{0} = \bfu_0$ and $\sobfu{1} = \bfu_0 + \eta \bfu_1$ in the electric case. In particular, $\sobfu{1}$ solves, due to linearity, the initial value problem
\begin{equation*}
\begin{aligned}
	\varepsilon_{\hom} \partial_{tt} \sobfu{1} - \Delta \sobfu{1} & = 0 \quad && \text{in } \R^d \times \R_+, \\
	\sobfu{1}(\cdot,0) & = \ve_0 \quad && \text{in } \R^d, \\
	\partial_t \sobfu{1}(\cdot,0) & = \ve_1 - \eta \chi_0 \varepsilon_{\hom}^{-1} \Delta \ve_0 \quad && \text{in } \R^d.
\end{aligned}
\end{equation*}
Considering $\sobfu{2}$, we see that it solves the homogenized equation
\begin{align} \label{eq:macro_u2}
	\varepsilon_{\hom} \partial_{tt} \sobfu{2} - \Delta \sobfu{2} & = \eta^2 \varepsilon_{\cor} \Delta^2 \bfu_0.
\end{align}
In particular, solving for $\sobfu{2}$ again requires solving for $\bfu_0$ first. Therefore, we aim for a $\mathcal{O}(\eta^3)$-approximation $\sbfu{2} = \sobfu{2} + \mathcal{O}(\eta^3)$, which solves an effective equation that does not depend directly on the quantities $\bfu_j$. We can rewrite \eqref{eq:macro_u2} using $\bfu_0 = \sobfu{2} - \eta \obfu_1 - \eta^2 \obfu_2$ and obtain
\begin{align} \label{eq11}
	\varepsilon_{\hom} \partial_{tt} \sobfu{2} - \Delta \sobfu{2} - \eta^2 \varepsilon_{\cor} \Delta^2 \sobfu{2} = - \eta^3 \varepsilon_{\cor} \Delta^2 \obfu_1 - \eta^4 \varepsilon_{\cor} \Delta^2 \obfu_2.
\end{align}
Neglecting the $\mathcal{O}(\eta^3)$-terms on the right hand side of \eqref{eq11} we can define $\sbfu{2}$ as the solution of the initial value problem
\begin{equation} \label{eq:pure_macro2_electric}
\begin{aligned}
	\varepsilon_{\hom} \partial_{tt} \sbfu{2} & - \Delta \sbfu{2} - \eta^2 \varepsilon_{\cor} \Delta^2 \sbfu{2} = 0 \quad && \text{in } \R^d \times \R_+, \\
	 \sbfu{2}(\cdot,0) & = \ve_0 - \eta^2 \varepsilon_{\hom}^{-1} \theta_0 \Delta \ve_0 \quad && \text{in } \R^d, \\
	\partial_t  \sbfu{2}(\cdot,0) & = \ve_1 - \eta \chi_0 \varepsilon_{\hom}^{-1} \Delta \ve_0 + \eta^2 \varepsilon_{\hom}^{-1} \theta_0 \Delta \ve_1\quad && \text{in } \R^d.
\end{aligned}
\end{equation}
We emphasize that $\sbfu{2}$ can be computed without pre-computation of any other quantity $\bfu_j$ and is a third order approximation of $\sobfu{2}$, i.e., $\sbfu{2} = \sobfu{2} + \mathcal{O}(\eta^3)$. In particular, $\sbfu{2} + \eta^2 \tbfu_2$ is a $\mathcal{O}(\eta^3)$-approximation of full wave solution $\bfu_\eta$.

\subsection{The electric field $\bfE$}

\begin{figure}[h]
\centering
\begin{minipage}{0.45\textwidth}
\centering
\includegraphics[scale=0.33]{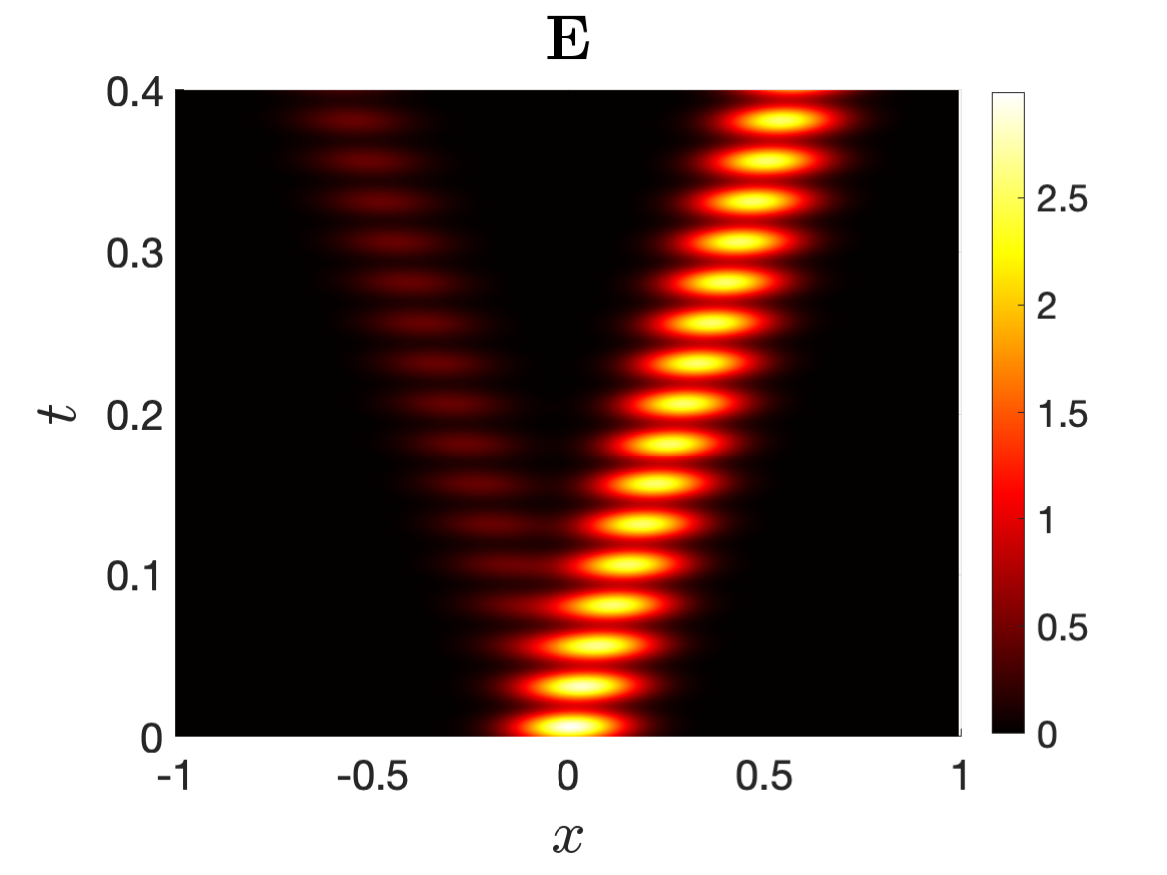}
\end{minipage}
\begin{minipage}{0.45\textwidth}
\centering
\includegraphics[scale=0.33]{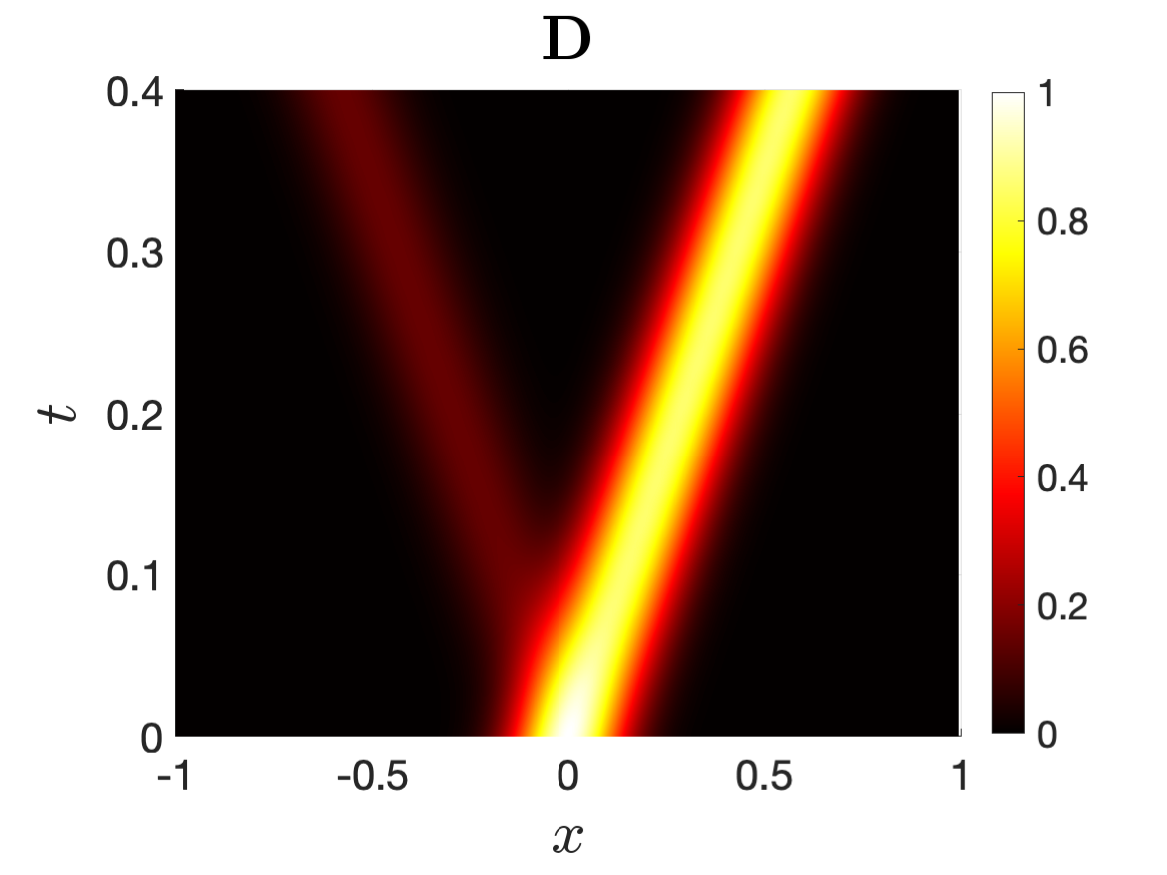}
\end{minipage}
\caption{Full-wave simulation of the electric fields $\bfE$ (left) and $\bfD$ (right) with permittivity $\varepsilon_\eta(t) = \big( 2 + \sin(2\pi t/\eta) \big)^{-1}$, $\eta = 0.025$ and initial values $\ve_0 = w(x,0)$, $\ve_1 = \partial_t w(x,0)$ where $w(x,t) =e^{-200(t - x)^2} \cos(\omega_0(t - x))$, $\omega_0 = 0.01$.} \label{fig:E_and_D}
\end{figure}

In the derivation of the model problems from Maxwell's equations in Section \ref{sec:setting} we also obtained an equation for the electric field $\bfE$, namely \eqref{eq:wave3}, for which the second time derivative in the wave equation acts on $\varepsilon \mathbf{E}$, and for which we could derive a similar effective behavior as in the previous sections. In particular, we obtain \eqref{eq:wave3} from \eqref{eq:wave2} by replacing $\bfD$ with $\varepsilon \bfE$ and multiplying with $\varepsilon^{-1}$. Since this connection is based on physical laws given by the constitutive relation \eqref{eq:material}, it should persist through the homogenization process. Homogenization approaches for \eqref{eq:wave3} reveal that $\bfD = \varepsilon \bfE$ is the quantity that can be homogenized and described by macroscopic quantities, but not $\bfE$. To illustrate this, we show in Figure \ref{fig:E_and_D} (left) an example of a $\bfE$ field solution from \eqref{eq:wave3}. The solution clearly shows microscopic oscillations on the length scale $\eta$ and therefore cannot be described by a homogenized $\eta$ independent quantity, while the quantity $\bfD = \varepsilon \bfE$ can, as discussed in the previous sections.

\section{The magnetic case} \label{subsec:hom:case1}

Next, we homogenize the magnetic field and consider the initial value problem \eqref{eq:case1}, where the electric permittivity occurs between the two time derivatives in the wave equation. Pretending that there is only one spatial dimension, we can exchange space and time, so that this case is similar to the homogenization of the (one-dimensional) wave equation with spatially multiscale coefficients (in divergence form), and we expect to obtain similar results. We follow the same strategy as in the electric case, plugging the asymptotic expansion into the PDE and initial conditions \eqref{eq:case1} and comparing powers in $\eta$. However, due to the different type of wave equations, we also get different problems for the effective solution and the higher order corrections compared to the electric case. In addition, three distinct cell problems on the unit interval $(0,1)$ arise describing the microscopic behavior of the first- and second-order correctors. To facilitate the reading of the subsequent derivation, we collect the cell problems and their respective roles in Table \ref{tab:cell-problems}. 

\begin{table}[h!]
    \centering
    \begin{tabular}{c c c}
        \toprule
        solution & cell-problem & role/determines \\
        \midrule
        $\chi$ & 
        $\partial_\tau(\varepsilon\partial_\tau\chi)=-\partial_\tau \varepsilon, \quad \displaystyle\int_0^1 \chi \, d\tau = 0$ & 
        \begin{minipage}[t]{0.4\textwidth}       
        \begin{itemize}
            \item microscopic components $\tbfu_1$ and $\tbfu_2$
            \item interface conditions for $\obfu_1$ and $\obfu_2$
        \end{itemize}        
        \end{minipage} \\
        \addlinespace
        $\xi$ & 
        $\partial_\tau (\varepsilon \partial_\tau \xi) = - \partial_\tau(\varepsilon\chi), \quad \displaystyle\int_0^1 \xi \, d\tau = 0$ & 
        \begin{minipage}[t]{0.4\textwidth}
        \begin{itemize}
            \item microscopic component $\tbfu_2$
            \item interface condition for $\obfu_2$
            \item $\xi(0) = -\theta_0$ with $\theta_0$ from \eqref{eq:theta_0}
        \end{itemize}  
        \end{minipage} \\
        \addlinespace
        $\zeta$ & 
        $\partial_\tau ( \varepsilon \partial_\tau \zeta) = - \partial_\tau(\varepsilon \xi) + \chi \varepsilon_{\hom}, \quad \displaystyle\int_0^1 \zeta \, d\tau = 0$ & 
        \begin{minipage}[t]{0.4\textwidth}
        \begin{itemize}
            \item only used as auxiliary problem
        \end{itemize}
        \end{minipage} \\
        \bottomrule
    \end{tabular}
    \caption{Overview of the cell-problems arising in the derivation of the magnetic homogenized equations.}
    \label{tab:cell-problems}
\end{table}

\subsection{Effective magnetic behavior} \label{sec:mag_eff}

Again we start with the derivation of a system of equations for the effective component $\bfu_0$ in \eqref{eq:asymptoticexpansion}. As in Section \ref{subsec:hom:case2}, we start considering the PDE  first and discuss the initial value afterwards. Plugging in the asymptotic expansion \eqref{eq:asymptoticexpansion} into the PDE in \eqref{eq:case1} yields
\begin{align*}
	0 & = \partial_t(\varepsilon_\eta\partial_t\bfu_\eta) - \Delta \bfu_\eta \\
	& =\frac{1}{\eta^2}\partial_\tau(\varepsilon\partial_\tau\bfu_0)+\frac{1}{\eta}\Bigl(\partial_t(\varepsilon\partial_\tau \bfu_0)+\partial_\tau(\varepsilon\partial_t\bfu_0)+\partial_\tau(\varepsilon\partial_\tau \bfu_1)\Bigr)\\
	& \quad + \sum_{j = 0}^{k-2} \eta^j \Bigl(\partial_t(\varepsilon\partial_t\bfu_j)+\partial_t(\varepsilon\partial_\tau \bfu_{j+1})+\partial_\tau(\varepsilon\partial_t\bfu_{j+1})+\partial_\tau(\varepsilon\partial_\tau\bfu_{j+2})- \Delta \bfu_j \Bigr) +\mathcal O(\eta^{k-1}).
\end{align*}
We compare powers of $\eta$ and first deduce $\partial_\tau(\varepsilon\partial_\tau\bfu_0)=0$ for the $\mathcal{O}(\eta^{-2})$-term. With the periodicity of $\bfu_0$ in $\tau$ this yields again $\tbfu_0 = 0$, i.e., $\bfu_0$ is independent of $\tau$. Considering the $\mathcal{O}(\eta^{-1})$-term we obtain
\begin{align} \label{eq12}
	\partial_\tau(\varepsilon\partial_t\bfu_0)+\partial_\tau(\varepsilon\partial_\tau \bfu_1)=0.
\end{align}
In contrast to the electric case, \eqref{eq12} does not imply $\bfu_1$ to be independent of $\tau$. Rather, we calculate (due to linearity) $\bfu_1$ in dependence on $\bfu_0$ and obtain
\begin{align*}
	\bfu_1(x,t,\tau)=\chi(\tau) \partial_t\bfu_0(x,t) + \overline{\bfu}_1(x,t).
\end{align*}
Here $\chi$ is the unique $1$-periodic solution of the first cell-problem, c.f. Table \ref{tab:cell-problems},
\begin{equation}\label{eq:case1-cellpb}
	\partial_\tau(\varepsilon\partial_\tau\chi)=-\partial_\tau \varepsilon, \quad \int_0^1 \chi \, d\tau = 0
\end{equation}
and $\overline{\bfu}_1$ still needs to be determined. In particular, $\chi$ from \eqref{eq:case1-cellpb} is explicitly given by
\begin{equation}\label{eq:w-explicit}
	\chi(\tau)=\chi_0 + \varepsilon_{\hom} \int_0^{\tau}\varepsilon^{-1} ds -\tau
\end{equation}
with $\chi_0$ from \eqref{chi0} such that $\int_0^1 \chi d\tau = 0$. Next, the $\mathcal{O}(\eta)$-term gives
\begin{equation} \label{eq:mag_eta0term}
	\partial_t(\varepsilon\partial_t\bfu_0)+\partial_t(\varepsilon\partial_\tau\bfu_1)+\partial_\tau(\varepsilon\partial_t\bfu_1) +\partial_\tau(\varepsilon\partial_\tau\bfu_2) -\Delta\bfu_0=0.
\end{equation} 
Integrating \eqref{eq:mag_eta0term} from $0$ to $1$ w.r.t. $\tau$ leads to
\begin{align} \label{eq13}
		\int_0^1 \partial_t(\varepsilon\partial_t\bfu_0)+\partial_t(\varepsilon\partial_\tau\bfu_1) -\Delta\bfu_0 \, d\tau = 0,
\end{align}
where we used the periodicity of $\varepsilon$,$\bfu_1$, and $\bfu_2$. Since $\bfu_0$ is independent of $\tau$ and $\partial_\tau \bfu_1 = (\partial_t\bfu_0) \partial_\tau \chi$, we deduce from \eqref{eq13}
\begin{equation*}
	\partial_t(\varepsilon_{\hom}\partial_t\bfu_0)-\Delta \bfu_0=0
\end{equation*}
with the homogenized coefficient
\begin{equation}\label{eq:case1:homcoeff}
\varepsilon_{\hom}\coloneqq \int_0^1\varepsilon(1+\partial_\tau \chi)\, d\tau.
\end{equation}
Inserting our explicit representation of $\chi$ from \eqref{eq:w-explicit} into \eqref{eq:case1:homcoeff}, we see that the expression from \eqref{eq:case1:homcoeff} coincides with the one from \eqref{eq:case1:harmav}, i.e., the homogenized coefficient is
\begin{equation*}
	\varepsilon_{\hom}=\Bigl(\int_0^1 \varepsilon^{-1} d\tau\Bigr)^{-1}.
\end{equation*}
Again, the homogenized coefficient is the harmonic average of the time-modulated permittivity over one time period. \\
We continue with the initial values. From the initial condition $\bfu_\eta(\cdot,0) = \bfv_0$ we immediately infer the first initial value $\bfu_0(\cdot,0) = \vm_0$. From the second initial condition $\varepsilon_\eta(0) \partial_t \bfu_\eta(\cdot,0) = \vm_1$ we find, after plugging in the asymptotic expansion, that $\varepsilon(0)(\partial_t \bfu_0(x,0) + \partial_\tau \bfu_1(x,0,0)) = \vm_1(x)$. Using $\partial_\tau \bfu_1(x,t,\tau) = \partial_\tau \chi(\tau) \partial_t \bfu_0(x,t) $ we obtain
\begin{align*}
	\varepsilon(0) (1 + \partial_\tau \chi(0)) \partial_t \bfu_0(x,0)= \vm_1(x). 
\end{align*}
By \eqref{eq:case1-cellpb} we conclude that $\varepsilon(1 + \partial_\tau \chi)$ is constant in $\tau$ and it is equal to $\varepsilon_{\hom}$ by \eqref{eq:case1:homcoeff}. This gives the second initial condition
\begin{align*}
	\varepsilon_{\hom} \partial_t \bfu_0(\cdot,0) = \vm_1.
\end{align*}
We summarize our findings regarding the effective solution $\bfu_0$:

\begin{Conclusion} \label{con:magnetic_hom}
Suppose the asymptotic expansion \eqref{eq:asymptoticexpansion} holds for the solution $\bfu_\eta$ of the initial value problem \eqref{eq:case1} up to order $k = 2$. Then $\tbfu_0 = 0$ in $\R^d \times \R_+$ and the effective solution $\bfu_0 = \obfu_0$ solves the homogenized initial value problem
\begin{equation} \label{eq:case1-homogenized}
\begin{aligned}
	\partial_t (\varepsilon_{\hom} \partial_t \bfu_0) - \Delta \bfu_0 & = 0 \quad && \text{in } \R^d \times \R_+, \\
	\bfu_0(\cdot,0) & = \vm_0 \quad && \text{in } \R^d, \\
	\varepsilon_{\hom} \partial_t \bfu_0(\cdot,0) & = \vm_1 \quad && \text{in } \R^d
\end{aligned}
\end{equation}
with the homogenized coefficient $\varepsilon_{\hom }$ from \eqref{eq:case1:harmav}.
\end{Conclusion}

\subsection{Magnetic first order correction} \label{sec:mag_first}

Let us now derive a full description of the first order correction $\bfu_1$. We already derived in the previous Section \ref{sec:mag_eff} that it is of the form
\begin{align} \label{eq:cor1_form}
	\bfu_1(x,t,\tau) = \chi(\tau) \partial_t \bfu_0(x,t) + \overline{\bfu}_1(x,t),
\end{align}
where $\chi$ is the solution of the cell-problem \eqref{eq:case1-cellpb} explicitly given by \eqref{eq:w-explicit}. Thus it remains to derive an equation for the time average $\overline{\bfu}_1$. For this purpose, we rewrite the $\mathcal{O}(1)$-term from \eqref{eq:mag_eta0term} and plug in \eqref{eq:cor1_form}, so that
\begin{align*}
	\partial_\tau (\varepsilon \partial_\tau \bfu_2) & = \Delta \bfu_0 - \partial_t(\varepsilon\partial_t \bfu_0) - \partial_t(\varepsilon \partial_\tau \bfu_1) - \partial_\tau (\varepsilon \partial_t \bfu_1) \\
	& = \Delta \bfu_0 - \varepsilon (1 + \partial_\tau \chi) \partial_{tt} \bfu_0 - \partial_\tau(\varepsilon \chi) \partial_{tt} \bfu_0 - (\partial_\tau \varepsilon )( \partial_t \overline{\bfu}_1) \\
	& = - \partial_\tau(\varepsilon\chi) \partial_{tt} \bfu_0 - (\partial_\tau \varepsilon )( \partial_t \overline{\bfu}_1).
\end{align*}
Here we used in the last step that $\bfu_0$ solves the homogenized equation \eqref{eq:case1-homogenized}. The equation is solved by
\begin{align*}
	\bfu_2(x,t,\tau) = \xi(\tau) \partial_{tt}\bfu_0(x,t) + \chi(\tau) \partial_t \overline{\bfu}_1(x,t)  + \overline{\bfu}_2(x,t),
\end{align*}
where $\xi$ is the unique $1$-periodic solution of the second cell-problem, c.f. Table \ref{tab:cell-problems},
\begin{align} \label{eq:cell-prob2}
	\partial_\tau (\varepsilon \partial_\tau \xi) = - \partial_\tau(\varepsilon\chi), \quad \int_0^1 \xi \, d\tau = 0.
\end{align}
The solution $\xi$ is simply given by 
\begin{align*}
	\xi(\tau) & = - \theta_0 -\int_0^\tau \chi \, ds,
\end{align*}
where we recall $\theta_0$ from \eqref{eq:theta_0}. Next we consider the $\mathcal{O}(\eta)$-term in the asymptotic expansion which yields
\begin{align} \label{eq20}
	\partial_\tau (\varepsilon \partial_\tau \bfu_3) + \partial_t (\varepsilon \partial_\tau \bfu_2) + \partial_\tau (\varepsilon \partial_t \bfu_2) + \partial_t (\varepsilon \partial_t \bfu_1) - \Delta \bfu_1= 0.
\end{align}
Integrating \eqref{eq20} over $(0,1)$ w.r.t. $\tau$ gives
\begin{align*}
	0 & = \int_0^1 \partial_t (\varepsilon \partial_{\tau} \bfu_2) + \partial_t (\varepsilon \partial_t \bfu_1) - \Delta \bfu_1 d\tau \\
	& = \int_0^1 \varepsilon (\chi + \partial_\tau \xi) \partial_{ttt} \bfu_0 + \varepsilon (1 + \partial_\tau \chi)  \partial_{tt} \overline{\bfu}_1 - \chi \partial_t \Delta \bfu_0 - \Delta \overline{\bfu}_1 \, d\tau.
\end{align*}
Since $\chi + \partial_\tau \xi \equiv 0$ and $\int_0^1 \chi d\tau = 0$, $\overline{\bfu}_1$ also solves the homogenized equation
\begin{align*}
	\partial_t ( \varepsilon_{\hom} \partial_{t} \overline{\bfu}_1) - \Delta \overline{\bfu}_1 =  0.
\end{align*}
The $\mathcal{O}(\eta)$-term in the first initial condition, i.e., $\bfu_\eta(\cdot,0) = \bfv_0$, implies $\bfu_1(x,0) = 0$ and therefore $\overline{\bfu}_1(\cdot,0) = -\chi_0 \partial_{t} \bfu_0(\cdot,0) = - \chi_0 \varepsilon_{\hom}^{-1} \vm_1$. For the second initial condition, $\varepsilon(0) \partial_t \bfu_\eta(\cdot,0) = \vm_1$, we obtain the $\mathcal{O}(\eta)$-term
\begin{align} \label{eq21}
	0 = \partial_t \bfu_1(x,0,\tau) + \partial_\tau \bfu_2(x,0,\tau).
\end{align}
Integration of \eqref{eq21} over $(0,1)$ w.r.t. $\tau$ yields due to the periodicity of $\bfu_2$
\begin{align*}
	0 = \int_0^1 \partial_t \bfu_1(x,0,\tau) d\tau = \int_0^1 \chi(\tau) d\tau \, \partial_{tt} \bfu_0(x,0) + \partial_t \overline{\bfu}_1(x,0) = \partial_t \overline{\bfu}_1(x,0). 
\end{align*}
We now summarize our findings for the first order correction $\bfu_1$:

\begin{Conclusion} \label{con:magnetic_cor1}
Suppose the asymptotic expansion \eqref{eq:asymptoticexpansion} holds for the solution $\bfu_\eta$ of the initial value problem \eqref{eq:case1} up to order $k = 3$. Then the first order correction $\bfu_1 = \tbfu_1 + \obfu_1$ is given by
\begin{equation*}
	\tbfu_1(x,t,\tau) = \chi(\tau) \partial_{t} \bfu_0(x,t)
\end{equation*}
with the solution $\chi$ of the cell-problem \eqref{eq:case1-cellpb}. Further, with the homogenized coefficient $\varepsilon_{\hom}$ from \eqref{eq:case1:harmav}, $\overline{\bfu}_1$ is given by the solution to the first order corrector problem
\begin{equation} \label{eq:case1-correction}
\begin{aligned}
	\partial_t (\varepsilon_{\hom} \partial_t \overline{\bfu}_1) - \Delta \overline{\bfu}_1 & = 0 \quad && \text{in } \R^d \times \R_+, \\
	\overline{\bfu}_1(\cdot,0) & = -\chi_0 \varepsilon_{\hom}^{-1} \vm_1 \quad && \text{in } \R^d, \\
	\varepsilon_{\hom} \partial_t \overline{\bfu}_1(\cdot,0) & = 0 \quad && \text{in } \R^d.
\end{aligned}
\end{equation}
\end{Conclusion}

\subsection{Magnetic higher order corrections}

Finally, we derive the system of equations for the second order correction $\bfu_2$. We have already shown in the previous section that $\bfu_2$ can be written as
\begin{align*}
	\bfu_2(x,t,\tau) = \widetilde{\bfu}_2(x,t,\tau) + \overline{\bfu}_2(x,t),
\end{align*}
where
\begin{align*}
	\widetilde{\bfu}_2(x,t,\tau) = \xi(\tau) \partial_{tt}\bfu_0(x,t) + \chi(\tau) \partial_{t} \overline{\bfu}_1(x,t).
\end{align*}
We consider the $\mathcal{O}(\eta)$-term in the asymptotic expansion, plug in the expressions for $\bfu_1$, $\bfu_2$, and use the fact that $\overline{\bfu}_1$ solves \eqref{eq:case1-correction} to obtain
\begin{align} \label{eq:mag_eta_term}
\begin{split}
	\partial_{\tau} (\varepsilon \partial_\tau \bfu_3) & = \Delta \bfu_1 - \partial_t (\varepsilon \partial_t \bfu_1) - \partial_t (\varepsilon \partial_\tau \bfu_2) - \partial_\tau (\varepsilon \partial_t \bfu_2) \\
	& = \Delta \overline{\bfu}_1 - \partial_t \varepsilon_{\hom} \partial_t \overline{\bfu}_1 + \chi \Delta \partial_t \bfu_0 - \partial_{\tau}(\varepsilon \xi) \partial_{ttt} \bfu_0 - \partial_\tau(\varepsilon \chi) \partial_{tt} \overline{\bfu}_1 - \partial_\tau \varepsilon \partial_t \overline{\bfu}_2 \\
	& = - (\partial_{\tau}(\varepsilon \xi) - \chi \varepsilon_{\hom}) \partial_{ttt} \bfu_0 - \partial_\tau(\varepsilon \chi) \partial_{tt} \overline{\bfu}_1 - \partial_\tau \varepsilon \partial_t \overline{\bfu}_2.
\end{split}
\end{align}
Here we replaced $\Delta \partial_t \bfu_0$ by $\varepsilon_{\hom} \partial_{ttt} \bfu_0$ in the last step. We can solve this for $\bfu_3$ and obtain
\begin{align*}
	\bfu_3(x,t,\tau) = \zeta(\tau) \partial_{ttt} \bfu_0(x,t) + \xi(\tau) \partial_{tt} \overline{\bfu}_1(x,t) + \chi(\tau) \partial_t \overline{\bfu}_2(x,t) + \overline{\bfu}_3(x,t),
\end{align*}
where $\zeta$ is the unique $1$-periodic solution of the third cell-problem, c.f. Table \ref{tab:cell-problems},
\begin{align*}
	\partial_\tau ( \varepsilon \partial_\tau \zeta) = - \partial_\tau(\varepsilon \xi) + \chi \varepsilon_{\hom}, \quad \int_0^1 \zeta \, d\tau = 0.
\end{align*}
Next we recall the $\mathcal{O}(\eta^2)$-term reading as
\begin{align} \label{eq23}
	\partial_\tau (\varepsilon \partial_\tau \bfu_4) + \partial_t (\varepsilon \partial_\tau \bfu_3) + \partial_\tau (\varepsilon \partial_t \bfu_3) + \partial_t (\varepsilon \partial_t \bfu_2) - \Delta \bfu_2= 0.
\end{align}
Integration of \eqref{eq23} over $(0,1)$ w.r.t. $\tau$ yields with our derived expression for $\bfu_3$ and $\bfu_2$
\begin{align*}
	0 & = \int_0^1 \partial_t(\varepsilon \partial_\tau \bfu_3) + \partial_t (\varepsilon \partial_t \bfu_2) - \Delta \bfu_2 d\tau \\
	& = \int_0^1 \varepsilon(\xi + \partial_\tau \zeta) \partial_t^4 \bfu_0 + \varepsilon(\chi + \partial_\tau \xi) \partial_t^3 \overline{\bfu}_1 + \varepsilon(1+ \partial_\tau \chi) \partial_{tt} \overline{\bfu}_2 - \Delta \overline{\bfu}_2 d\tau \\
	& = \int_0^1 \varepsilon_{\hom}^{-2} \varepsilon(\xi + \partial_\tau \zeta) \Delta^2 \bfu_0 + \varepsilon(\chi + \partial_\tau \xi) \partial_t^3 \overline{\bfu}_1 + \varepsilon(1+ \partial_\tau \chi) \partial_{tt} \overline{\bfu}_2 - \Delta \overline{\bfu}_2 d\tau,
\end{align*}
where in the last step we replaced $\partial_t^4 \bfu_0$ by $\varepsilon_{\hom}^{-2} \Delta^2 \bfu_0$ if $\bfu_0$ is a smooth solution of \eqref{eq:case1-homogenized}. Next, we make the important observation that, by straightforward calculations,
\begin{align*}
	\varepsilon_{\hom}^{-2} \int_0^1 \varepsilon(\xi + \partial_\tau \zeta) d\tau = - \int_0^1 \varepsilon^{-1} \theta \, d\tau = - \varepsilon_{\cor}.
\end{align*}
Therefore, with $\chi + \partial_\tau \xi \equiv 0$, $\overline{\bfu}_2$ also solves the homogenized equation with the corrected right-hand side
\begin{align*}
	\partial_t (\varepsilon_{\hom} \partial_t \overline{\bfu}_2) - \Delta \overline{\bfu}_2 = \varepsilon_{\cor} \Delta^2 \bfu_0 
\end{align*}
and with the correction coefficient from \eqref{electric_correction_coefficient}. For the initial conditions for $\bfu_2$, respectively $\overline{\bfu}_2$, we observe from $\bfu_\eta(\cdot,0) = \vm_0$ that $\bfu_2(x,0,0) = 0$ and therefore 
\begin{align*}
	\overline{\bfu}_2(x,0) = \theta_0 \partial_{tt} \bfu_0(x,0) - \chi_0 \partial_{t} \overline{\bfu}_1(x,0) = \theta_0\partial_{tt} \bfu_0(x,0) = \varepsilon_{\hom}^{-1} \theta_0 \Delta \vm_0,
\end{align*}
since $\partial_t \overline{\bfu}_1(x,0) = 0$ and where we assume that $\bfu_0$ is a smooth solution of \eqref{eq:case1-correction2}. For the second initial condition, $\varepsilon(0) \partial_t \bfu_\eta(\cdot,0) = \vm_1$, we now obtain with the $\mathcal{O}(\eta^2)$-term in the asymptotic expansion
\begin{align*}
	0 = \partial_t \bfu_2(x,0,0) + \partial_\tau \bfu_3(x,0,0).
\end{align*}
Recalling the $\mathcal{O}(\eta)$-term \eqref{eq:mag_eta_term} in the asymptotic expansion, we solve for $\partial_\tau \bfu_3$ and obtain
\begin{align} \label{eq24}
\begin{split}
	\varepsilon(\tau) \partial_\tau \bfu_3(x,t,\tau) - \varepsilon(0) \partial_\tau \bfu_3(x,t,0) & = \Big(\varepsilon(0) \xi(0) - \varepsilon(\tau) \xi(\tau)+ \varepsilon_{\hom} \int_0^\tau \chi \, ds \Big) \partial_{ttt} \bfu_0(x,t) \\
	& \quad + ( \varepsilon(0) \chi_0 - \varepsilon(\tau) \chi(\tau)) \partial_{tt} \obfu_1(x,t) + (\varepsilon(0) - \varepsilon(\tau)) \partial_t \obfu_2(x,t). 
\end{split}
\end{align}
Dividing \eqref{eq24} by $\varepsilon(\tau)$ and by $\varepsilon(0)$ yields, after integration over $(0,1)$ w.r.t. $\tau$,
\begin{align} \label{eq25}
\begin{split}
	-  \partial_\tau  \bfu_3(x,t,0) & =  -\theta_0 \partial_{ttt} \bfu_0(x,t) + \frac{\varepsilon_{\hom}^2}{\varepsilon(0)} \int_0^1 \varepsilon^{-1} \int_0^\tau \chi \, ds \, \partial_{ttt} \bfu_0(x,t) \\
	& \quad + \chi_0 \partial_{tt} \obfu_1(x,t) + \partial_t \obfu_2(x,t) - \frac{\varepsilon_{\hom}}{\varepsilon(0)} \partial_t \obfu_2(x,t).
\end{split}
\end{align}
With $\partial_t \obfu_2(x,0) =  - \partial_\tau \bfu_3(x,0,0) - \partial_t \tbfu_2(x,0,0)$ and $\partial_t \tbfu_2(x,0,0) = -\theta_0 \partial_{ttt} \bfu_0(x,0) + \chi_0 \partial_{tt} \obfu_1(x,0)$ equation \eqref{eq25} implies
\begin{align*}
	\varepsilon_{\hom} \partial_t \obfu_2(x,0) & = \varepsilon_{\hom}^2 \Big( \int_0^1 \varepsilon^{-1} \int_0^\tau \chi \, ds d\tau \Big) \partial_{ttt} \bfu_0(x,0) = \Big( \int_0^1 \varepsilon^{-1} \int_0^\tau \chi \, ds d\tau \Big) \Delta \vm_1.
\end{align*}
We summarize our findings for the second order correction $\bfu_2$:

\begin{Conclusion} \label{con:magnetic_cor2}
Suppose the asymptotic expansion \eqref{eq:asymptoticexpansion} holds for the solution $\bfu_\eta$ of the initial value problem \eqref{eq:case1} up to order $k = 4$. Then the second order correction $\bfu_2 = \tbfu_2 + \obfu_2$ is given by
\begin{equation*}
	\tbfu_2(x,t,\tau) = \xi(\tau) \partial_{tt} \bfu_0(x,t) +  \chi(\tau) \partial_{t} \overline{\bfu}_1(x,t) 
\end{equation*}
with the solutions $\xi,\chi$ of the cell-problems \eqref{eq:cell-prob2} and \eqref{eq:case1-cellpb}. Furthermore, with the homogenized coefficient $\varepsilon_{\hom}$ from \eqref{eq:case1:harmav} and the correction coefficient $\varepsilon_{\cor}$ from \eqref{electric_correction_coefficient}, $\overline{\bfu}_2$ is given by the solution to the second order corrector problem
\begin{equation} \label{eq:case1-correction2}
\begin{aligned}
	\partial_t (\varepsilon_{\hom} \partial_t \overline{\bfu}_2) - \Delta \overline{\bfu}_2 & = - \varepsilon_{\cor} \Delta^2 \bfu_0  \quad && \text{in } \R^d \times \R_+, \\
	\overline{\bfu}_2(\cdot,0) & = \varepsilon_{\hom}^{-1} \theta_0 \Delta \vm_0 \quad && \text{in } \R^d, \\
	\varepsilon_{\hom} \partial_t \overline{\bfu}_2(\cdot,0) & = \Big( \int_0^1 \varepsilon^{-1} \int_0^\tau \chi \, ds d\tau \Big) \Delta \vm_1  \quad && \text{in } \R^d.
\end{aligned}
\end{equation}
\end{Conclusion}

\subsection{Pure macroscopic description of the magnetic field}
\label{sec:pure_marco_magnetic}

As in the electric case, let us consider the macroscopic quantities given by \eqref{macro_quantities} for the magnetic field. Again we have $\sobfu{0} = \bfu_0$ and the initial value problem for $\sobfu{1} = \bfu_0 + \eta \bfu_1$ now reads as
\begin{equation*}
\begin{aligned}
	\partial_t (\varepsilon_{\hom} \partial_t \sobfu{1}) - \Delta \sobfu{1} & = 0 \quad && \text{in } \R^d \times \R_+, \\
	\sobfu{1}(\cdot,0) & = \vm_0 - \eta \chi_0 \varepsilon_{\hom}^{-1} \vm_1 \quad && \text{in } \R^d, \\
	\varepsilon_{\hom}  \partial_t \sobfu{1}(\cdot,0) & = \vm_1 \quad && \text{in } \R^d.
\end{aligned}
\end{equation*}
In addition, we derive an initial value problem for a third order approximation $\sbfu{2}$ of $\sobfu{2}$, i.e., $\sbfu{2} = \sobfu{2} + \mathcal{O}(\eta^3)$ given by
\begin{equation}\label{eq:pure_macro2_magnetic}
\begin{aligned}
	\partial_t (\varepsilon_{\hom} \partial_t \sbfu{2}) & - \Delta \sbfu{2} - \eta^2 \varepsilon_{\cor} \Delta^2 \sbfu{2} = 0 \quad && \text{in } \R^d \times \R_+, \\
	 \sbfu{2}(\cdot,0) & = \vm_0 - \eta \chi_0 \varepsilon_{\hom}^{-1} \vm_1 + \eta^2 \varepsilon_{\hom}^{-1} \theta_0 \Delta \vm_0 \quad && \text{in } \R^d, \\
	\varepsilon_{\hom}  \partial_t  \sbfu{2}(\cdot,0) & = \vm_1 + \eta^2 \Big( \int_0^1 \varepsilon^{-1} \int_0^\tau \chi \, ds d\tau \Big) \Delta \vm_1 \quad && \text{in } \R^d.
\end{aligned}
\end{equation}

We close this section with a discussion of our findings, focusing in particular on the comparison between the electric and magnetic cases. At the effective level, both fields satisfy the homogenized wave equation with the effective coefficient $\varepsilon_{\hom}$ and the respective pure initial values $\ve_0, \ve_1$ for the electric field, and $\vm_0, \vm_1$ for the magnetic field. The first notable differences appear in the first-order corrections: in the electric case, no microscopic dynamics are observed, whereas in the magnetic case, microscopic dynamics arise through the solution $\chi$ of the first cell problem \eqref{eq:case1-cellpb}. At the macroscopic level, both first-order corrections again satisfy the homogenized wave equation with coefficient $\varepsilon_{\hom}$, but now with modified initial values induced by the interface conditions (see Conclusions \ref{con:electric_cor1} and \ref{con:magnetic_cor1}). Finally, for the second-order corrections in both fields, microscopic dynamics are induced via solutions of higher-order cell problems. At the macroscopic level, we obtain, up to a sign, the same homogenized PDEs for the second-order components, now featuring higher-order (nonlocal) source terms and appropriately calculated initial values.

\section{Homogenized Maxwell's equations} \label{sec:Maxwell}

In this section we revise the system of Maxwell's equations \eqref{eq:Maxwell} in view of our observations from the previous sections. In particular, we derive the homogenized Maxwell systems for the effective behavior of an electromagnetic wave propagating through a time-varying metamaterial. The homogenized Maxwell systems for the first order corrections as well as the equations for nonlocal second order corrections are revised, revealing local and nonlocal material laws. \\  
Recall the Maxwell system \eqref{eq:Maxwell} which can be formulated with the constitutive relations \eqref{eq:material} with permittivity $\varepsilon_\eta$ as a system for the pair $(\bfD,\bfB)$ reading as
\begin{equation} \label{eq:Maxwell2}
	\begin{aligned}
	\curl\bfD+ \varepsilon_\eta \partial_t \bfB & = 0, \quad & \Div\bfD=0, \\
	\curl \bfB-\partial_t \bfD & =0, \quad & \Div\bfB=0.
	\end{aligned}
\end{equation}
This system is equipped with the initial conditions
\begin{align*}
	\bfD(x,0) = \ve_0, \quad \bfB(x,0) = \vm_0.
\end{align*}
It is now clear from our derivation in Section \ref{sec:setting} that $\bfD$ and $\bfB$ solve the initial wave-type problems \eqref{eq:case2} and \eqref{eq:case1}, respectively. Thus, we can apply our observation from Section \ref{subsec:hom:case2} and Section \ref{subsec:hom:case1} to the associated electric and magnetic field. For this purpose, we assume that the fields $\bfD$ and $\bfB$ can be formally expanded in $\eta$ as
\begin{align*}
	\bfD(x,t) & = \bfD_0(x,t,\tau) + \eta \bfD_1(x,t,\tau) + \eta^2 \bfD_2(x,t,\tau) + \cdots, \\
	\bfB(x,t) & = \bfB_0(x,t,\tau) + \eta \bfB_1(x,t,\tau) + \eta^2 \bfB_2(x,t,\tau) + \cdots,
\end{align*}
where $\tau = t/\eta$, $\bfD_j(x,t,\cdot)$, $\bfB_j(x,t,\cdot)$ are $1$-periodic and can be decomposed into macroscopic and microscopic components according to
\begin{align*}
	\bfD_j(x,t,\tau) = \obfD_j(x,t) + \tbfD_j(x,t,\tau), \quad \bfB_j(x,t,\tau) = \obfB_j(x,t) + \tbfB_j(x,t,\tau).
\end{align*}
Here $\tbfD_j$ and $\tbfB_j$ have zeromean over $\tau$; cf. \eqref{zero-mean}. From \eqref{eq:Maxwell2} we immediately conclude $\Div \bfD_j = 0$ and $\Div \bfB_j = 0$ for all $j \in \N$ and we can further assume w.l.o.g that $\Div \obfD_j = 0$ and $\Div \obfB_j = 0$ for all $j \in \N$. In view of the interface conditions of Morgenthaler \cite{Mor1958} (cf. \eqref{eq:interface1} and \eqref{eq:interface2})  we set $\ve_1 = \curl \vm_0$ and $\vm_1 = - \curl \ve_0$. In this setting, the fields $\bfD_j$ and $\bfB_j$ for $j = 0,1,2$ and its macro-microscopic components, respectively, are given through Conclusion \ref{con:electric_hom}-\ref{con:magnetic_cor2}.

\begin{description}
\item[Homogenized equations ($j = 0$):] Conclusion \ref{con:electric_hom} and Conclusion \ref{con:magnetic_hom} imply that the microscopic components $\tbfD_0, \tbfB_0$ vanish, i.e., $\tbfD_0 = \tbfB_0 = 0$. Further, the macroscopic components $\bfD_0, \bfB_0$ solve the initial value problems \eqref{eq:case2-homogenized} and \eqref{eq:case1-homogenized}, respectively. These can be reformulated as a homogeneous Maxwell system such that the pair $(\bfD_0, \bfB_0)$ solves
\begin{equation*}
	\begin{aligned}
	\curl\bfD_0+ \varepsilon_{\hom} \partial_t \bfB_0 & = 0, \quad & \Div\bfD_0=0, \\
	\curl \bfB_0-\partial_t \bfD_0 & =0, \quad & \Div\bfB_0=0.
	\end{aligned}
\end{equation*}
The initial conditions for this system are given by $\bfD_0(x,0) = \ve_0$ and $\bfB_0(x,0) = \vm_0$.

\item[First order correction ($j = 1$):] We can proceed by formulating the system for the first order corrections $\bfD_1$ and $ \bfB_1$. From Conclusion \ref{con:electric_cor1} and Conclusion \ref{con:magnetic_cor1} we infer that $\tbfD_1 = 0$ and
\begin{align*}
\tbfB_1(x,t,\tau) = -\varepsilon_{\hom}^{-1} \chi(\tau)   \curl \bfD_0(x,t)
\end{align*}
with $\chi$ from \eqref{eq:case1-cellpb}. The macroscopic components $\obfD_1$ and $\obfB_1$ now solve \eqref{eq:case2-firstcorrectionsystem} and \eqref{eq:case1-correction}, respectively. This can be reformulated as homogeneous Maxwell system for $(\bfD_1, \obfB_1)$ reading as
\begin{equation*}
	\begin{aligned}
		\curl\bfD_1+ \varepsilon_{\hom} \partial_t \obfB_1&= 0, \quad && \Div\bfD_1=0,\\
		\curl \obfB_1-\partial_t \bfD_1&=0, \quad && \Div\obfB_1=0\\
	\end{aligned}
\end{equation*}
with the initial conditions $\bfD_1(x,0) = 0$ and $\obfB_1(x,0) = - \chi_0 \varepsilon_{\hom}^{-1} \vm_1$.

\item[Second order correction ($j = 2$):] Next we treat the second order corrections and the nonlocal effects. We obtain from Conclusion \ref{con:electric_cor2} and Conclusion \ref{con:magnetic_cor2}
\begin{align*}
	\tbfD_2(x,t,\tau) & = -\varepsilon_{\hom}^{-1} \theta(\tau) \curl \curl \bfD_0(x,t), \\
	\tbfB_2(x,t,\tau) & = - \varepsilon_{\hom}^{-1} \xi(\tau) \curl \curl \bfB_0(x,t)  - \varepsilon_{\hom}^{-1} \chi(\tau) \curl \bfD_1(x,t)
\end{align*}
with $\theta$ from \eqref{electric_cell_problem}, $\xi$ from \eqref{eq:cell-prob2}, and $\chi$ from \eqref{eq:case1-cellpb}. Furthermore, the macroscopic components $\obfD_2$ and $\obfD_1$ solve \eqref{eq:case2-correction2} and \eqref{eq:case1-correction2}, respectively, but  clearly they do not satisfy a type of homogeneous Maxwell system. Instead, we consider the concatenated macroscopic components
\begin{align*}
	\sobfD{2}(x,t) & = \bfD_0(x,t) + \eta \bfD_{1}(x,t) + \eta^2 \obfD_2(x,t), \\
	\sobfB{2}(x,t) & = \bfB_0(x,t) + \eta \obfB_{1}(x,t) + \eta^2 \obfB_2(x,t),
\end{align*}
which are again divergence-free, i.e., $\Div \sobfD{2} = 0$ and  $\Div \sobfB{2} = 0$.
Then $\sobfD{2}$ solves \eqref{eq:macro_u2} so that
\begin{align} \label{eq36}
	\varepsilon_{\hom} \partial_{tt} \sobfD{2} - \Delta \sobfD{2} = \eta^2 \varepsilon_\cor \Delta^2 \bfD_0 = \eta^2 \varepsilon_{\hom} \varepsilon_\cor \Delta \partial_{tt} \bfD_0.
\end{align}
Replacing $\bfD_0$ by $\sobfD{2} - \eta \bfD_1 - \eta^2 \obfD_2$ in \eqref{eq36}, we obtain
\begin{align} \label{eq37}
	(\varepsilon_{\hom} - \varepsilon_{\hom} \varepsilon_\cor \Delta) \partial_{tt} \sobfD{2} - \Delta \sobfD{2} = \eta^3 \varepsilon_{\hom} \varepsilon_\cor \Delta \partial_{tt} \bfD_1 + \eta^4 \varepsilon_{\hom} \varepsilon_\cor \Delta \partial_{tt} \obfD_2.
\end{align}
Dropping the $\mathcal{O}(\eta^3)$ source terms on the right hand side of \eqref{eq37}, we define, in analogy of Section \ref{sec:pure_marco_electric}, the macroscopic electric quantity $\sbfD{2}$ as the solution of
\begin{align*}
	(\varepsilon_{\hom} - \varepsilon_{\hom} \varepsilon_\cor \Delta) \partial_{tt} \sbfD{2} - \Delta \sbfD{2} = 0
\end{align*}
with the initial conditions
\begin{align*}
	\sbfD{2}(x,0) & =\ve_0 - \eta^2 \varepsilon_{\hom}^{-1} \theta_0 \Delta \ve_0, \\
	\partial_t \sbfD{2}(x,0) & = \ve_1 - \eta \chi_0 \varepsilon_{\hom}^{-1} \Delta \ve_0 + \eta^2 \varepsilon_{\hom}^{-1} \theta_0 \Delta \ve_1.
\end{align*}
Then $\sbfD{2}$ is an $\mathcal{O}(\eta^3)$-approximation of $\sobfD{2}$, i.e., $\sbfD{2} = \sobfD{2} + \mathcal{O}(\eta^3)$. Similarly, we define the macroscopic magnetic quantity $\sbfB{2}$ as the solution of
\begin{align*}
	(\varepsilon_{\hom} - \varepsilon_{\hom} \varepsilon_\cor \Delta) \partial_{tt} \sbfB{2} - \Delta \sbfB{2} = 0
\end{align*}
with the initial conditions
\begin{align*}
	\sbfB{2}(x,0) & =\vm_0 - \eta \chi_0 \vm_1 + \eta^2 \varepsilon_{\hom}^{-1} \theta_0 \Delta \vm_0, \\
	\varepsilon_{\hom} \partial_t \sbfB{2}(x,0) & = \vm_1 + \eta^2 \Big( \int_0^1 \varepsilon^{-1} \int_0^\tau \chi \, ds d\tau \Big) \Delta \vm_1.
\end{align*}
This is an $\mathcal{O}(\eta^3)$-approximation of $\sobfB{2}$, i.e., $\sbfB{2} = \sobfB{2} + \mathcal{O}(\eta^3)$. Note that by construction
\begin{align*}
	\sbfD{2} + \eta^2 \tbfD_2 = \bfD + \mathcal{O}(\eta^3) \quad \text{and} \quad \sbfB{2} + \eta \tbfB_1 + \eta^2 \tbfB_2 = \bfB + \mathcal{O}(\eta^3).
\end{align*}
Now we obtain that the pair $(\sbfD{2}, \sbfB{2})$ is a solution of the Maxwell system
\begin{equation*}
	\begin{aligned}
		\curl\sbfD{2}+ (\varepsilon_{\hom} + \eta^2 \varepsilon_{\hom} \varepsilon_{\cor} \curl \curl ) \partial_t \sbfB{2}&= 0, \quad && \Div\sbfD{2}=0,\\
		\curl \sbfB{2}-\partial_t \sbfD{2}&=0 , \quad && \Div\sbfB{2}=0,
	\end{aligned}
\end{equation*}
with the initial conditions
\begin{align*}
		\sbfD{2}(x,0) & =\ve_0 - \eta^2 \varepsilon_{\hom}^{-1} \theta_0 \Delta \ve_0, \\
		\sbfB{2}(x,0) & =\vm_0 - \eta \chi_0 \vm_1 + \eta^2 \varepsilon_{\hom}^{-1} \theta_0 \Delta \vm_0.
\end{align*}
This leads to the nonlocal constitutive relation
\begin{equation*}
	\sbfD{2} = (\varepsilon_{\hom} + \eta^2 \varepsilon_{\hom} \varepsilon_{\cor} \curl \curl ) \mathbf{E}^{(2)}
\end{equation*}
between the macroscopic electric displacement $\sbfD{2}$ and the associated macroscopic electric field $\mathbf{E}^{(2)}$. In the limit $\eta \rightarrow 0$, this gives $\sbfD{2} = \varepsilon_{\hom} \mathbf{E}^{(2)}$, which is the same local constitutive relation as that of the effective fields. However, for $\eta > 0$ there is an additional contribution of $\eta^2 \varepsilon_{\hom} \varepsilon_{\cor} \curl \curl \mathbf{E}^{(2)}$. This additional contribution is nonlocal, as it involves derivatives of the electric field. For monochromatic plane waves, this nonlocal material law produces artificial magnetism in the material by altering the effective magnetic permeability obtained from the dispersion relation; cf. \cite{RCG22,Torrent20}. As demonstrated in the numerical simulations reported in \cite{RCG22}, incorporating the nonlocal term (artificial magnetism) into the model of light propagation can significantly improve the accuracy of the approximation in certain parameter regimes, compared to a purely local description. In the next section, we present analogous numerical observations in particular in the regime with $\omega_0 \ll \omega_m$, where we recall that $\omega_m = \eta^{-1}$ and $\omega_0$ denotes the spatial modulation frequency. For a quantitative analysis of the relevant parameter regimes, we refer to the experimental study in \cite{RCG22}.
\end{description}

\section{Numerical experiments}
\label{sec:numerics}

In this section, we verify our homogenization results in numerical simulations.
We aim to (i) illustrate the derived homogenized solutions for the electric and magnetic case and to (ii) verify the expected approximation orders in terms of $\eta$. The implementation of the experiments is available as a MATLAB code on \url{https://github.com/cdoeding/PropagationTimeVaryingMedia}. \\

\paragraph*{Model setting} 
As a model problem, we consider a transverse electromagnetic wave that propagates through a time-varying metamaterial, and we simulate its electric displacement and magnetic field independently. This allows us to describe the time evolution of the fields through the equations \eqref{eq:case2} and \eqref{eq:case1} in one dimension. Note that our homogenization results do not depend on the spatial dimension and, thus, the consideration of the one-dimensional case is not restrictive. We model the situation that the time oscillations of the material are instantly switched on at time $t = 0$. For negative times $t < 0$ we assume that the material behaves as vacuum ($\varepsilon(t) = 1$) so that the wave propagates through the medium at the speed of light $c = 1$. In this case the electric displacement or magnetic field, respectively, can be described by a Gaussian wave packet of the form 
\begin{equation}\label{eq:wavepacket}
	\bfw(x,t) = \exp\Bigl(-\frac{(t- x )^2}{2T_0^2}\Bigr)\cos(\omega_0(t- x )), \quad x \in \R, \, t \le 0,
\end{equation}
with $T_0 > 0$ and $\omega_0 \ge 0$. We choose $T_0 = 0.1$ and $\omega_0 = 0.01$ in the subsequent experiments. Note that $\omega_0$ is the carrier frequency and determines the spatial oscillations of the wave that we need to keep small in comparison to the modulation frequency $\omega_m = \eta^{-1}$ so that the homogenization errors are small in view of the asymptotic expansion \eqref{eq:asymptoticexpansion}. \\
For $t \ge 0$ the permittivity of the material is given by the blueprint
\begin{align} \label{blueprint_num}
	\varepsilon(\tau) = \big( 2 + \sin(2\pi \tau) \big)^{-1}
\end{align}
and we vary the time modulations through different values of $\eta$ satisfying $\omega_0 \ll \eta^{-1}$. The homogenized coefficient from \eqref{eq:case1:harmav} and correction coefficient \eqref{electric_correction_coefficient} can be calculated explicitly and are given by
\begin{align*} 
	\varepsilon_{\hom} = \frac{1}{2}, \quad \varepsilon_{\cor} = - \frac{1}{16 \pi^2}.
\end{align*}
Furthermore, the solution of the cell-problems \eqref{electric_cell_problem}, \eqref{eq:case1-cellpb}, and \eqref{eq:cell-prob2} can be solved explicitly. However, we note that in any particular application the homogenized coefficients $\varepsilon_{\hom}$ and $\varepsilon_\cor$ as well as the solutions to the cell-problems can be solved numerically using suitable quadratures if an explicit representation is not available. \\
For the initial values in our model we use that the interface conditions of Morgenthaler \cite{Mor1958} require the fields to be continuous at $t = 0$ leading to the first initial condition for the electric field $\ve_0 = \bfw(x,0)$ and for the magnetic field $\vm_0 = \bfw(x,0)$, respectively. A precise recap of the model and in particular its simplification from three to one dimension reveals that the second initial conditions for the wave-type problems \eqref{eq:case2} and \eqref{eq:case1} are given by $\ve_1 = \partial_t \bfw(x,0)$ and $\vm_1 = \partial_t \bfw(x,0)$, respectively. Indeed the latter two conditions coincide with the conditions \eqref{eq:interface1} and \eqref{eq:interface2}.

\paragraph*{Numerical discretization} 

What remains to be solved numerically are the time-averaged solutions $\obfu_j$ given by the initial-value problems stated in Conclusion \ref{con:electric_hom}-\ref{con:magnetic_cor2}. In order to solve the corresponding wave equations, we first restrict the initial-value problem to the bounded domain $\Omega = (-1,1)$, impose periodic boundary conditions, and solve on a finite time interval $[0,T]$ with $T = 0.4$. Note that the initial conditions and system parameters are chosen such that the propagating wave never reaches the boundary in the time interval $[0,T]$, so that the artificial boundary conditions do not induce any additional effects. In space, we discretize using a Fourier spectral method with $N = 256$ (for visualization) or $N = 64$ (for error computation) degrees of freedom, which allows a simple discretization of $\Delta$ and $\Delta^2$ in Fourier space and so that the discretization errors are negligible. The Fourier spectral method is combined with the two-stage Gauss-Legendre IRK of order four with variable time step size $\tau = 2^{-8} T$ for visualization and $\tau = 2^{-13}T$ for the error computation. Again, the high order of the time integrator allows us to neglect any errors induced by the time discretization, so that only the errors induced by the homogenization process become dominant in our experiments. The actual full wave solution $\bfu_\eta$ is solved with the same method and parameter sets and serves as a reference solution. \\

\paragraph*{Discussion of the results} 

\begin{figure}[h]
\centering
\begin{minipage}{0.325\textwidth}
\centering
\includegraphics[scale=0.25]{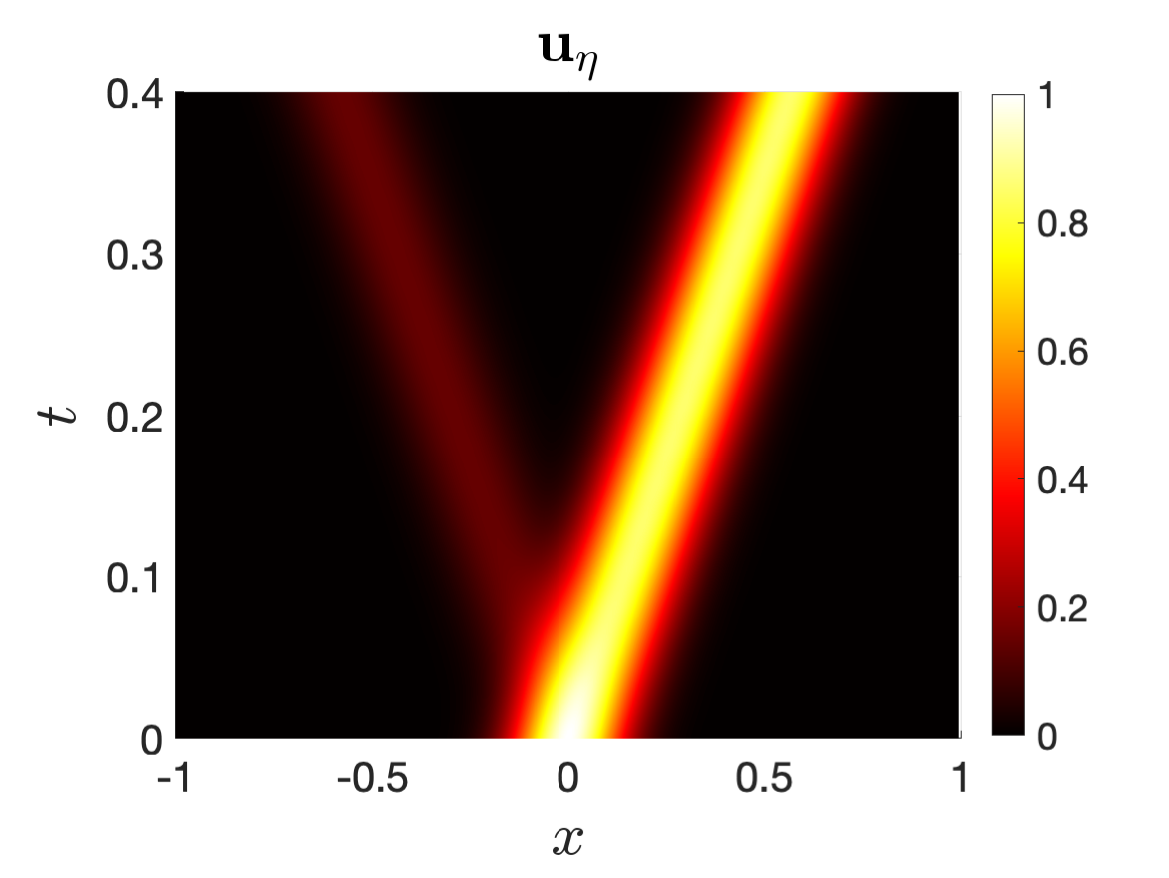}
\end{minipage}
\begin{minipage}{0.325\textwidth}
\centering
\includegraphics[scale=0.25]{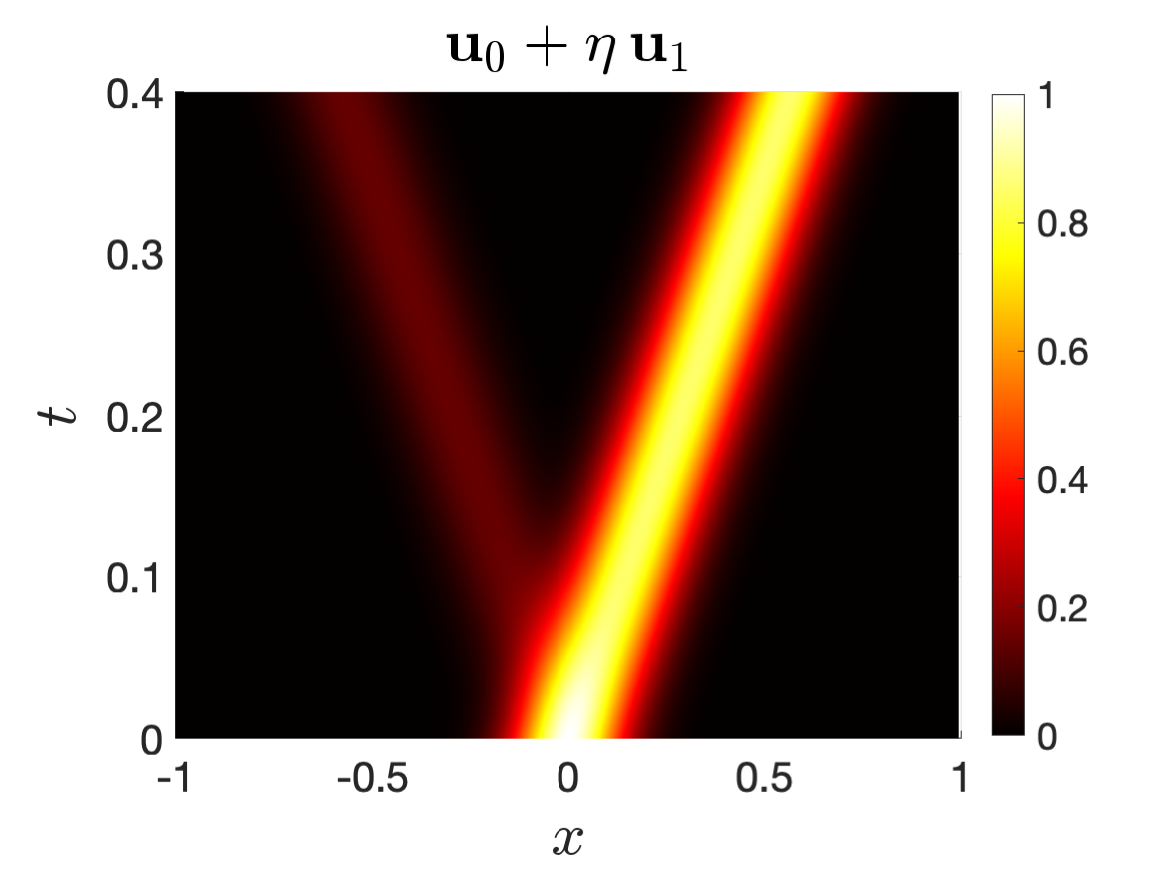}
\end{minipage}
\begin{minipage}{0.325\textwidth}
\centering
\includegraphics[scale=0.25]{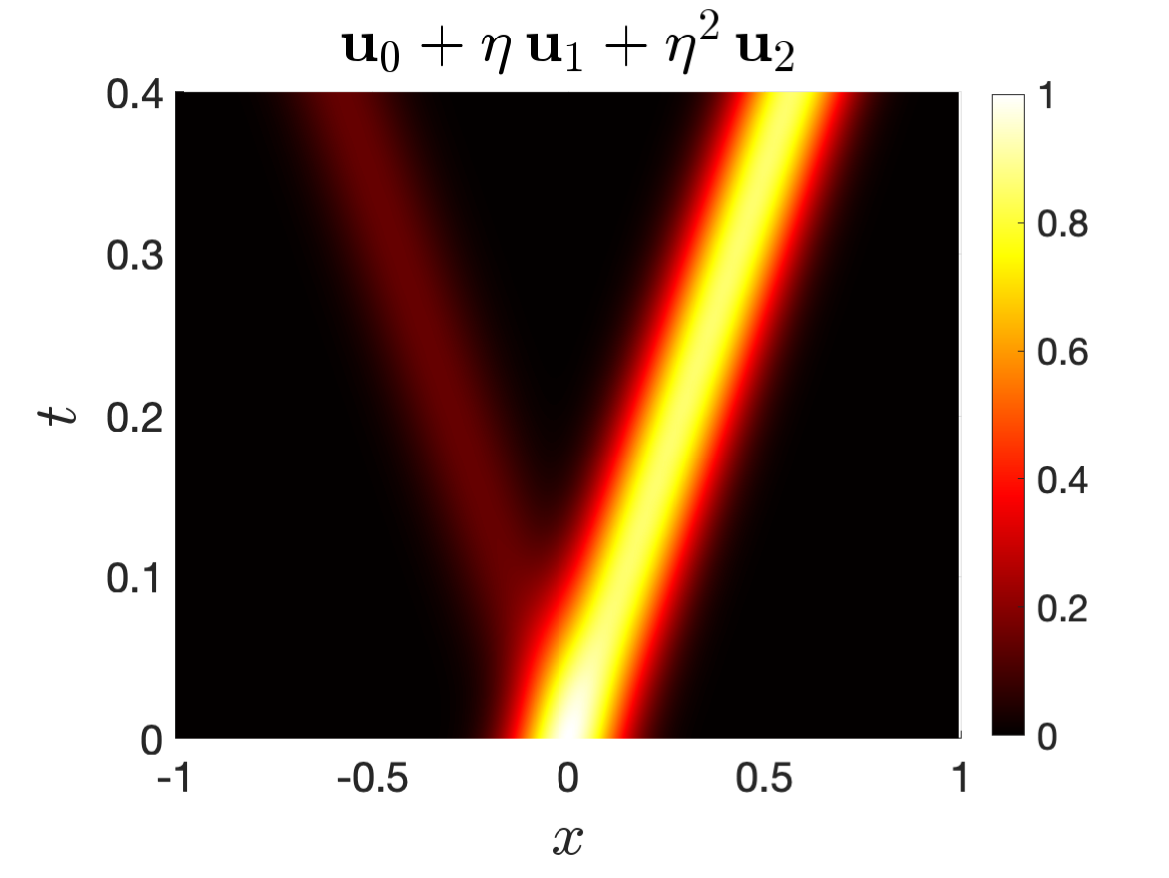}
\end{minipage} \\
\begin{minipage}{0.325\textwidth}
\centering
\includegraphics[scale=0.25]{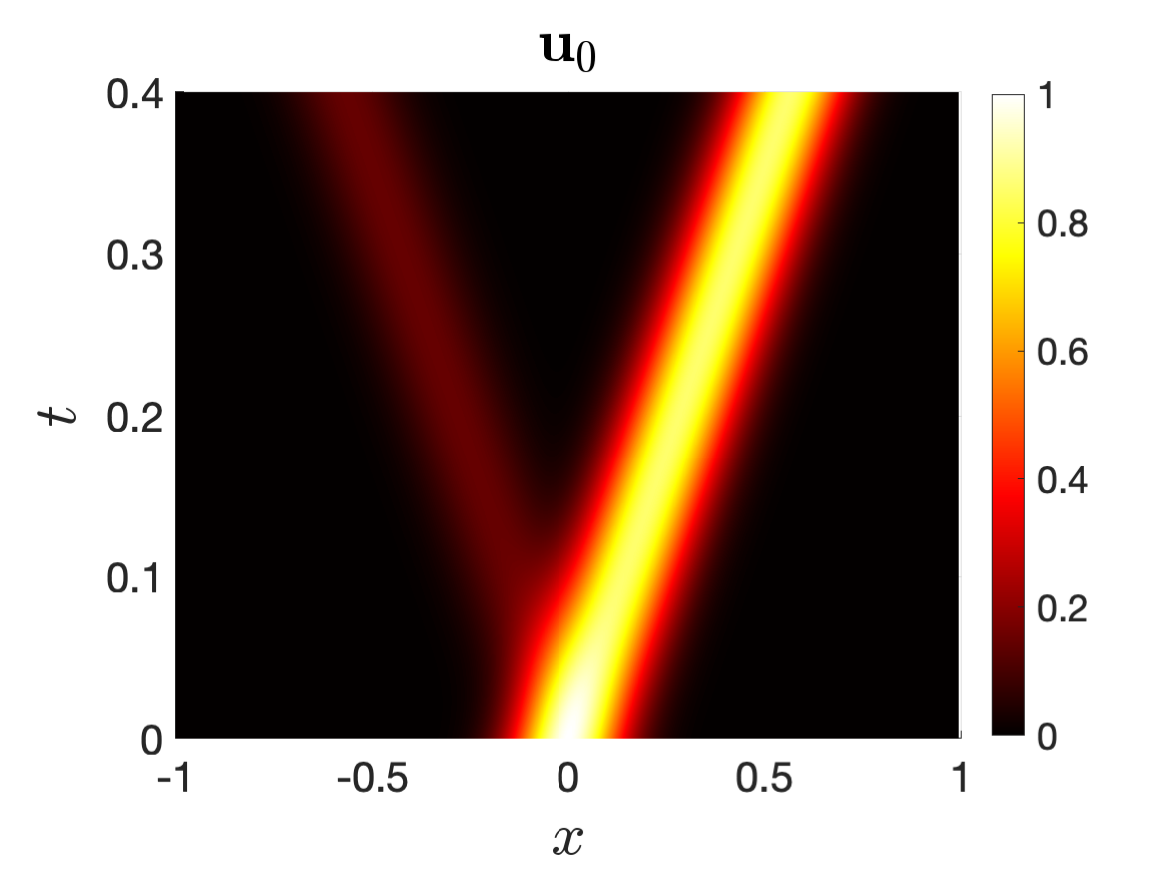}
\end{minipage}
\begin{minipage}{0.325\textwidth}
\centering
\includegraphics[scale=0.25]{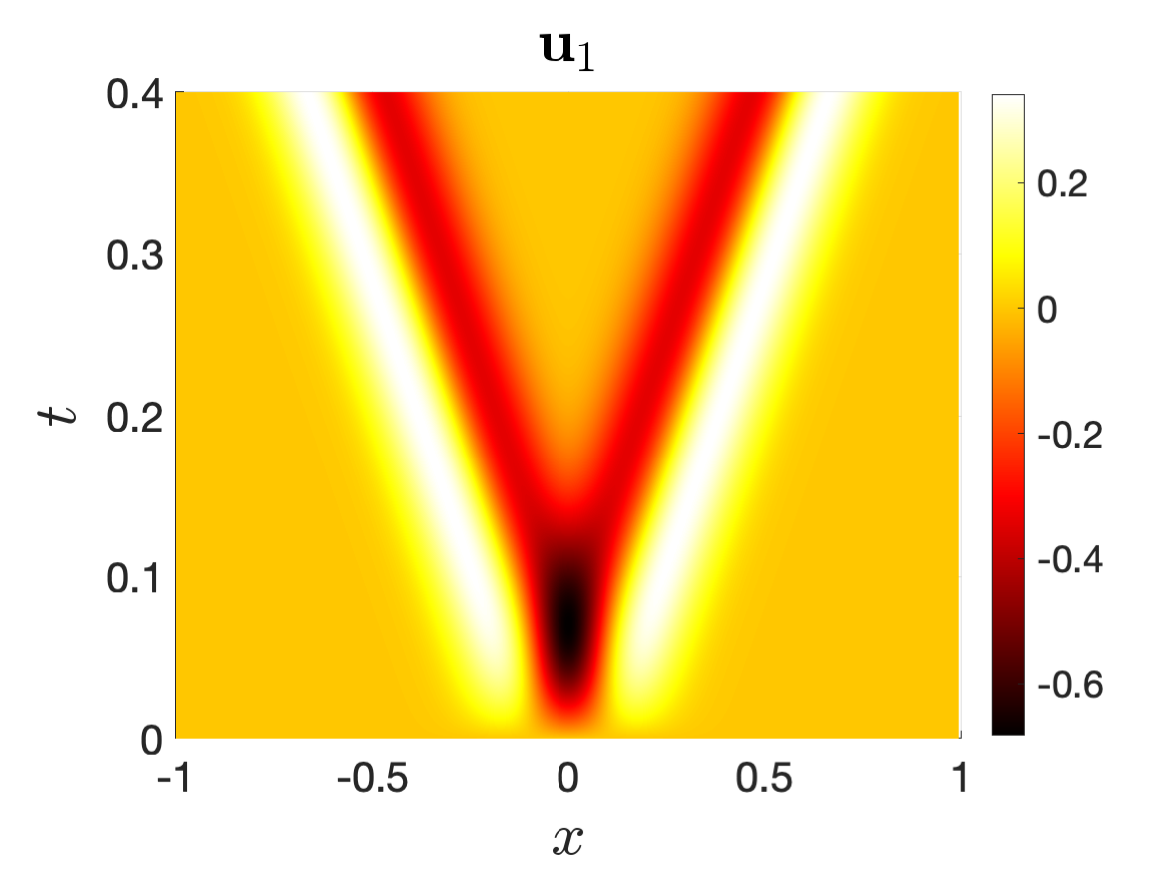}
\end{minipage}
\begin{minipage}{0.325\textwidth}
\centering
\includegraphics[scale=0.25]{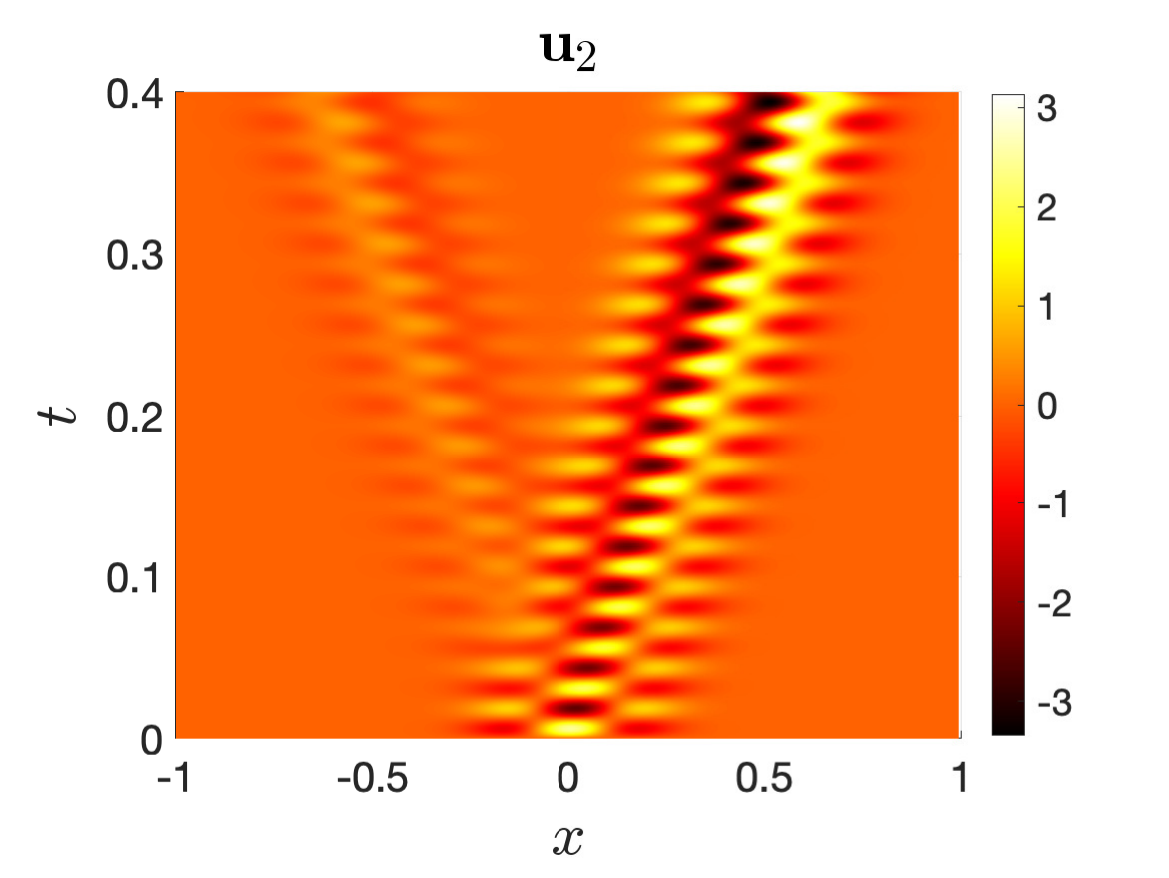}
\end{minipage}
\caption{Full wave simulation and effective behavior of the electric field for the fine-scale parameter $\eta = 2^{-4}T$.} \label{numerics:electric}
\end{figure}

\begin{figure}[h]
\centering
\begin{minipage}{0.325\textwidth}
\centering
\includegraphics[scale=0.26]{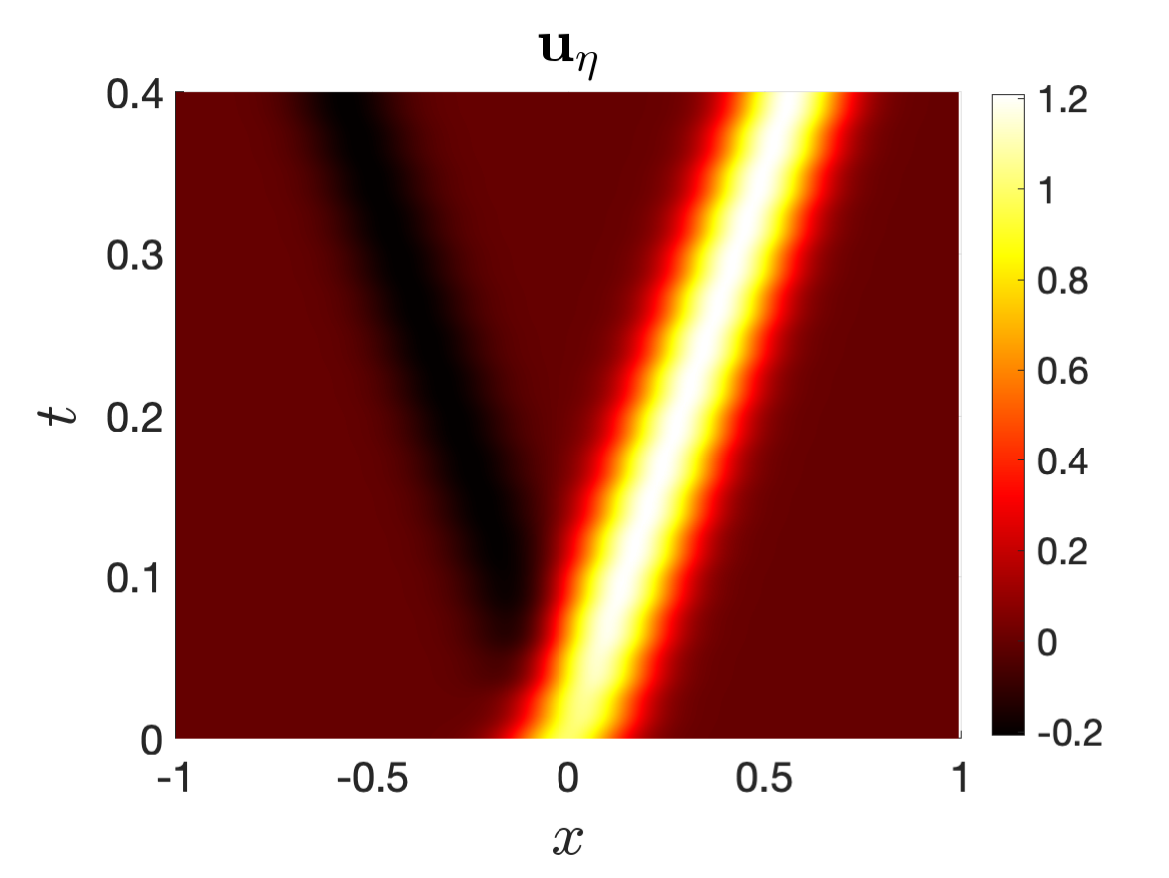}
\end{minipage}
\begin{minipage}{0.325\textwidth}
\centering
\includegraphics[scale=0.26]{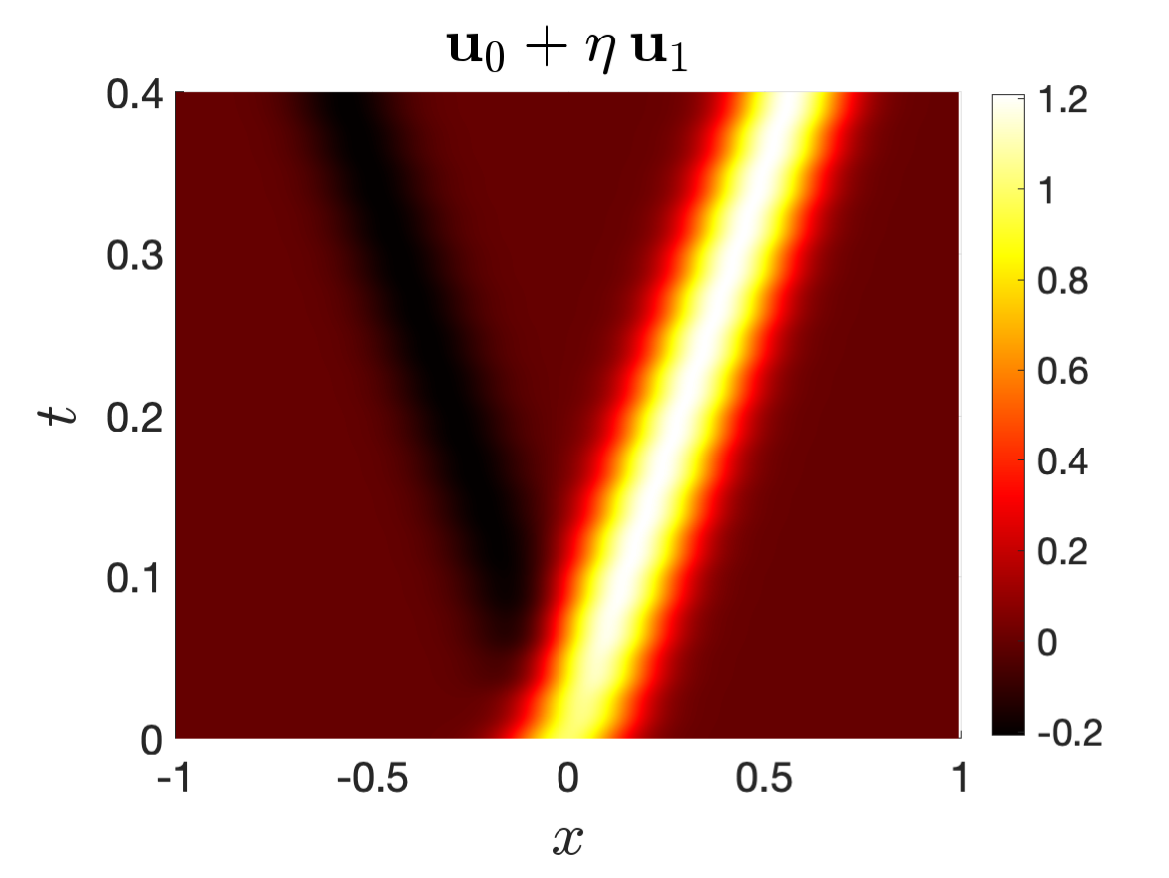}
\end{minipage}
\begin{minipage}{0.325\textwidth}
\centering
\includegraphics[scale=0.26]{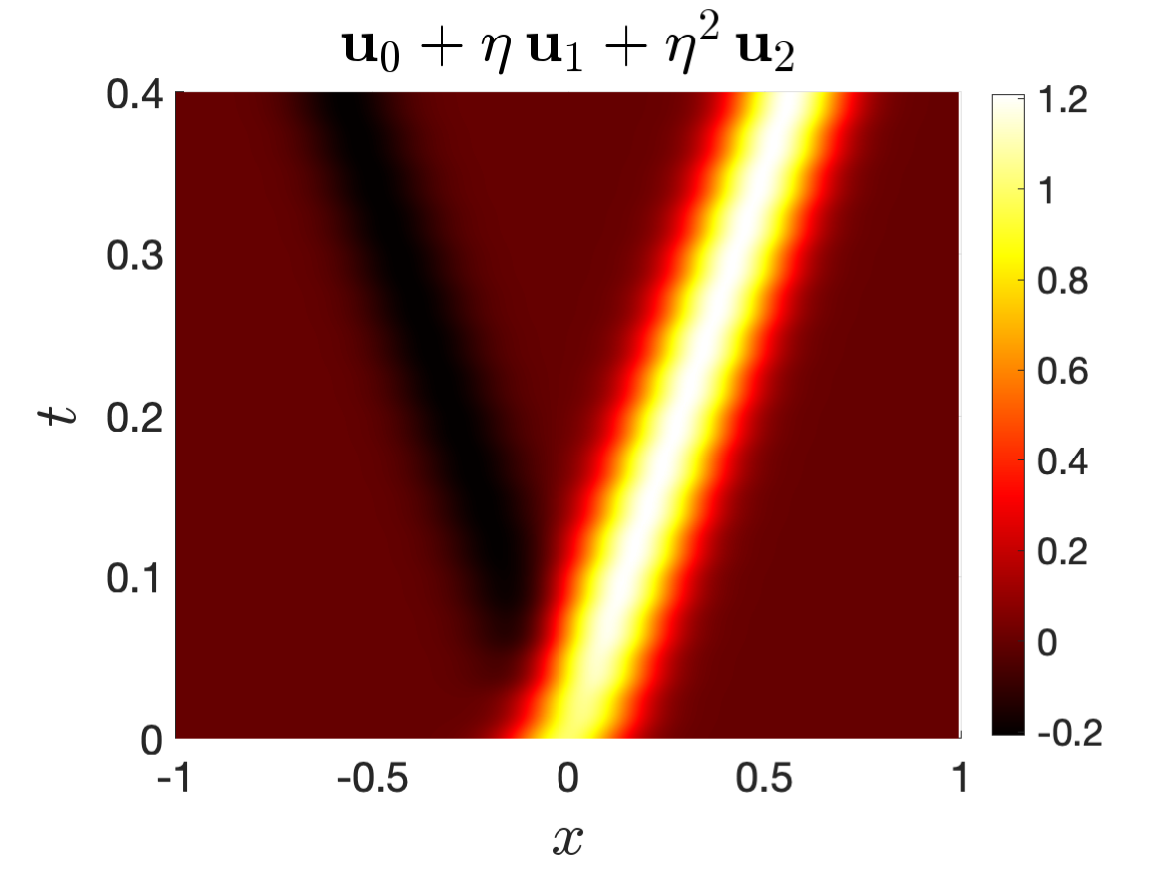}
\end{minipage} \\
\begin{minipage}{0.325\textwidth}
\centering
\includegraphics[scale=0.26]{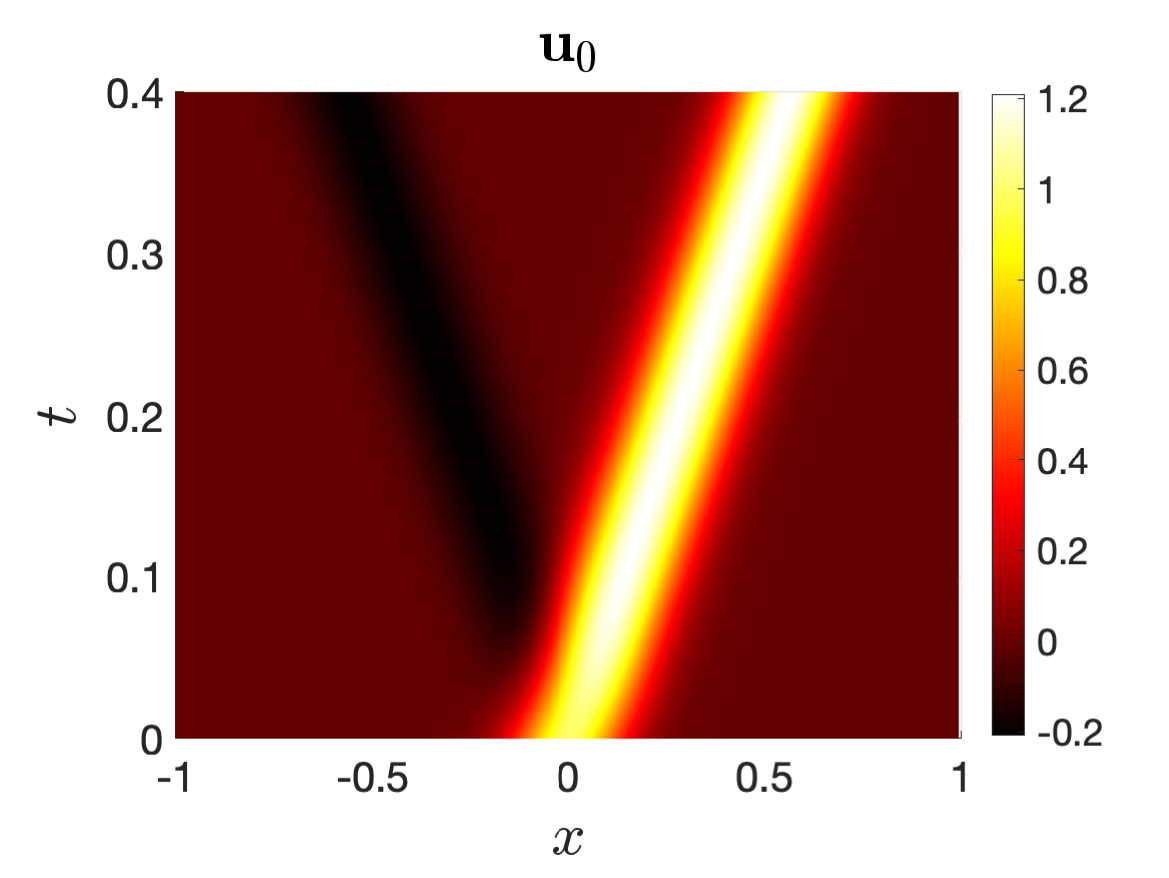}
\end{minipage}
\begin{minipage}{0.325\textwidth}
\centering
\includegraphics[scale=0.26]{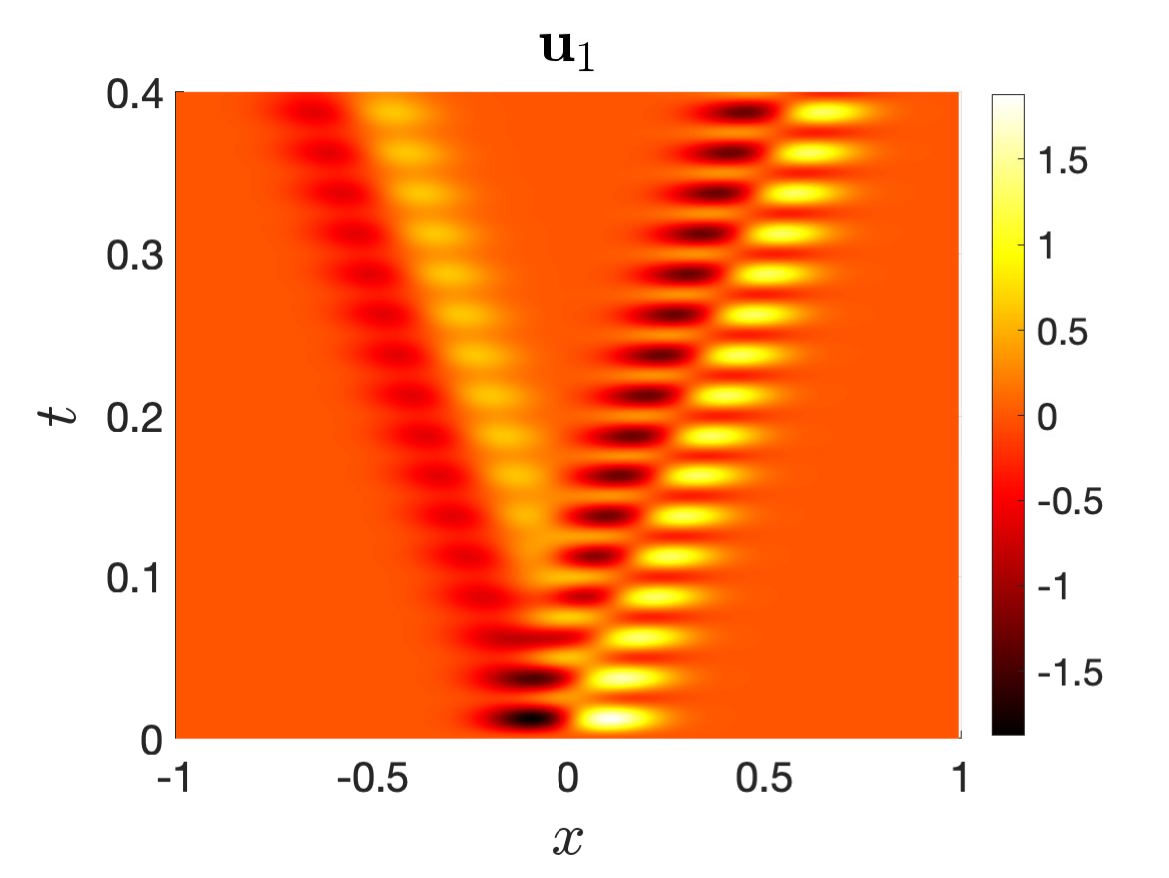}
\end{minipage}
\begin{minipage}{0.325\textwidth}
\centering
\includegraphics[scale=0.26]{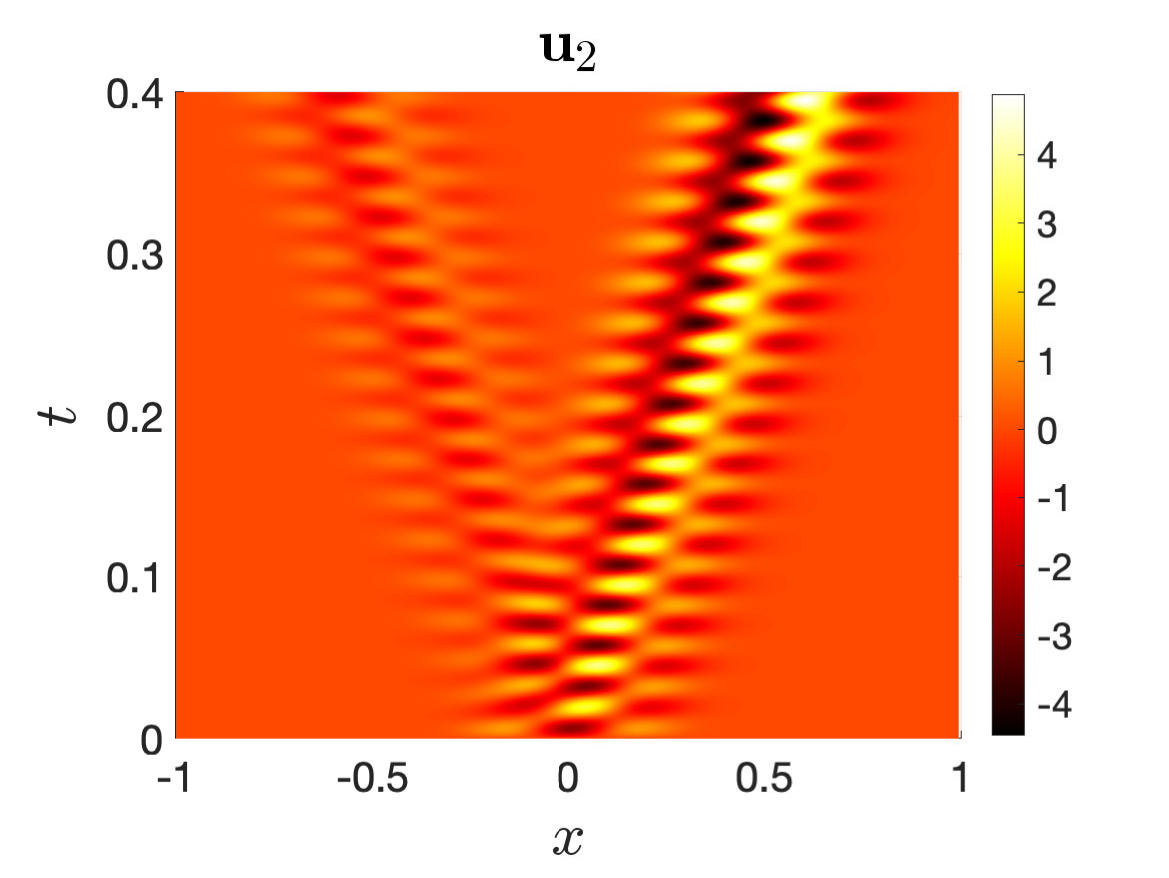}
\end{minipage}
\caption{Full wave simulation and effective behavior of the magnetic field for the fine-scale parameter $\eta = 2^{-4}T$.} \label{numerics:magnetic}
\end{figure}

We first consider the electric displacement and compute the full wave problem \eqref{eq:case2} for the choice $\eta = 2^{-4} T$. Here we use $N = 64$ Fourier modes, step size $\tau = 2^{-8}T$ and plot the full wave solution $\bfu_\eta$ in Figure \ref{numerics:electric} (top left). As expected, we observe the wave propagating to the right, accompanied by a another reflected wave traveling to the left. This second component can be identified as a temporal reflection phenomenon, which is well-known to arise at the temporal interface $t = 0$, see e.g. \cite{Galiffi22}. More specifically, it corresponds to a phase-conjugated component with negative frequency, generated when the metamaterial switches from a constant permittivity for $t < 0$ to a time-varying permittivity for $t> 0$. The temporal interface conditions thus split the wave into a forward-propagating (time-refracted) part and a backward-propagating (phase-conjugated) part. While no fine-scale oscillations are visually apparent, the propagation speed is clearly governed by the temporal modulation of the permittivity. \\
The effective solution $\bfu_0$ as given in Conclusion \ref{con:electric_hom} is shown in Figure \ref{numerics:electric} (bottom left). As expected, the effective solution captures the macroscopic behavior of to the full wave solution $\bfu_\eta$ and resolves the effective propagation of the wave through the metamaterial. The first and second order homogenized solutions $\bfu_1$ and $\bfu_2$ as given in Conclusion \ref{con:electric_cor1} and Conclusion \ref{con:electric_cor2} are also shown in Figure \ref{numerics:electric}. We clearly see that the first order correction $\bfu_1$ does not contain any microscopic information but gives a clear contribution to the transmitted and reflected wave. Only the second order correction $\bfu_2$ contains microscopic information in accordance with Conclusion \ref{con:electric_cor2}. However, visually at this point we see no difference between the corrected solutions $\bfu_0 + \eta \bfu_1$ and $\bfu_0 + \eta \bfu_1 + \eta^2 \bfu_2$ to the effective solution $\bfu_0$ and the full wave simulation $\bfu_\eta$. However, the effective wave propagation through the time varying metamaterial is captured by all homogenized solutions. \\
Next we consider the magnetic field of the electromagnetic wave described by the initial-value problem \eqref{eq:case1}, which we solve with the same set of parameters. The solution $\bfu_\eta$ is shown in Figure \ref{numerics:magnetic} and again we see a transmitted wave propagating to the right and a reflected wave propagating to the left. Unlike the electric case, the reflected wave now has negative values and, furthermore, shows a microscopic behavior through fine oscillations. These fine oscillations cannot be resolved by the effective solution $\bfu_0$ from Conclusion \ref{con:magnetic_hom}, as shown in Figure \ref{numerics:magnetic} (bottom left). Nevertheless, $\bfu_0$ captures again the effective behavior of the wave. In contrast, the first-order correction $\bfu_1$ from Conclusion \ref{con:magnetic_cor1} resolves the fine oscillation, leading to no visual difference between the corrected effective solution $\bfu_0 + \eta \bfu_1$ and the full wave solution $\bfu_\eta$. The second order correction shown in the Figure \ref{con:magnetic_cor2} (bottom right) adds further micro- and macroscopic information to the solution, but is visually indistinguishable from the first order corrected effective solution. \\
Summarizing, our theoretical results lead to visually good approximations of the full wave solution in both the electric and the magnetic case, and the effective and corrected equations resolve the propagation of the electromagnetic wave through the time-varying metamaterial. \\

\begin{figure}[h]
\centering
\begin{minipage}{0.49\textwidth}
\centering
\includegraphics[scale=0.4]{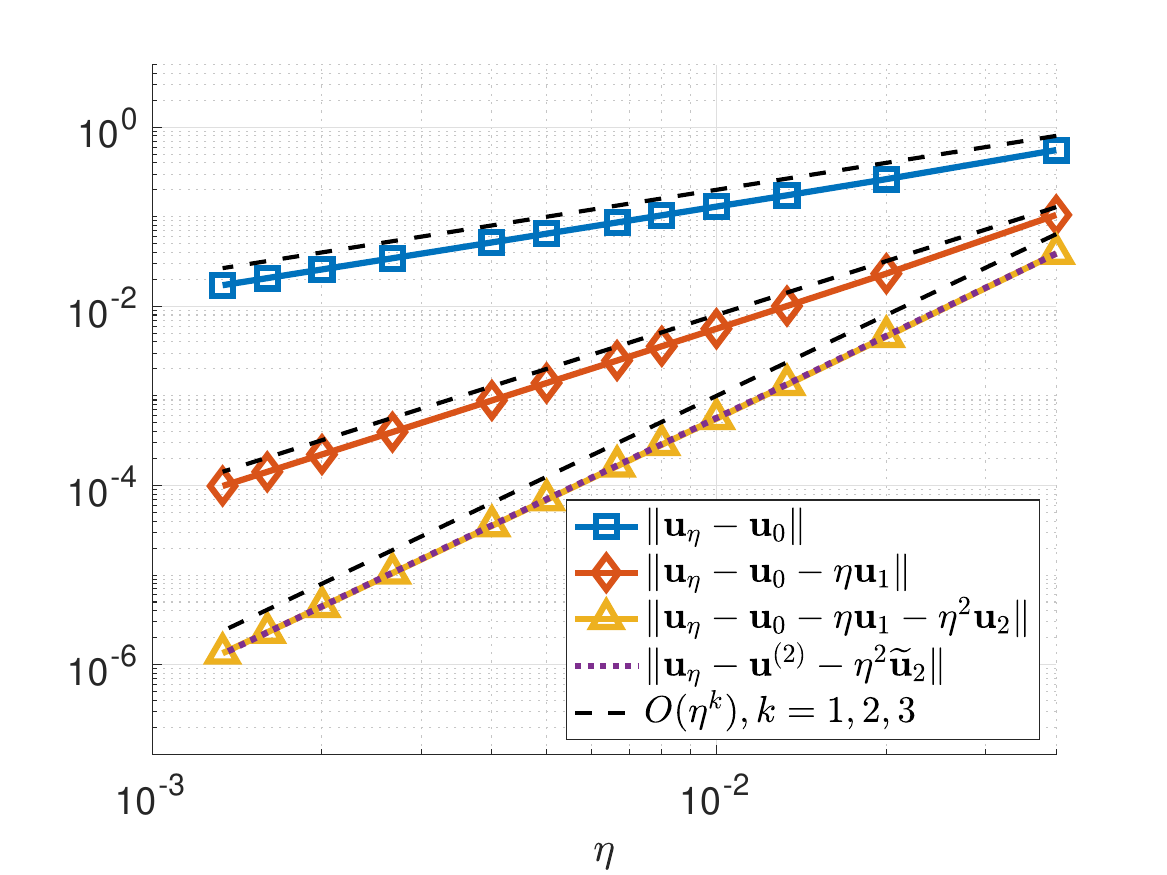}
\end{minipage}
\begin{minipage}{0.49\textwidth}
\centering
\includegraphics[scale=0.4]{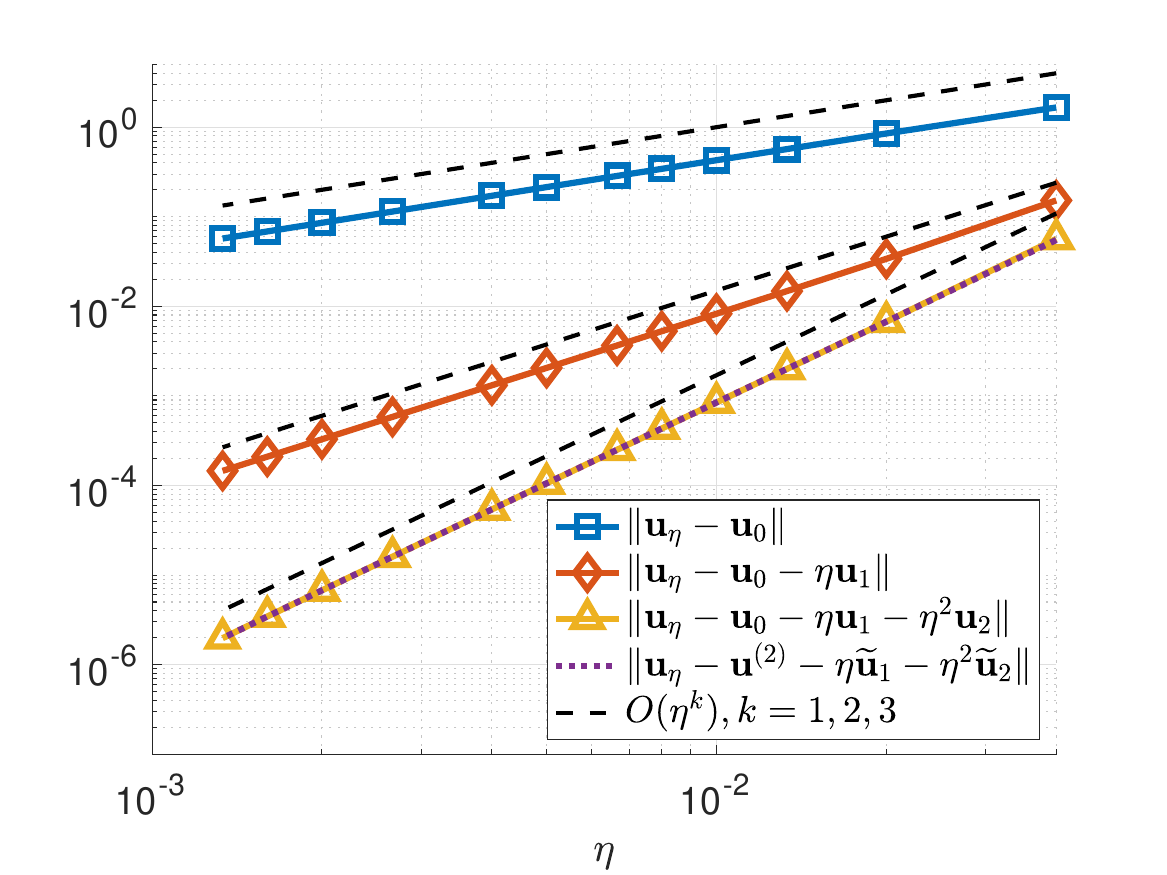}
\end{minipage}
\caption{Homegenization errors in $\eta$ for the electric field (left) and in the magnetic field (right).} \label{error_plot}
\end{figure}

To quantify the approximation properties of the homogenized solutions and to validate the derived equations in Section \ref{subsec:hom:case2} and Section \ref{subsec:hom:case1}, we study the homogenization error 
\begin{align*}
	\| \bfu_\eta - \bfu_0 \|_{L^2(0,T; L^2)}, 
\end{align*} 
for a set of fine-scale parameters $\eta = T/\ell$, $\ell = 10, 20, 30, 40, 50, 60, 80, 100, 150, 200, 250, 300$. Similarly, we consider the homogenization errors for the first and second order homogenized solutions $\bfu_0 + \eta \bfu_1$ and $\bfu_0 + \eta \bfu_1 + \eta^2 \bfu_2$. For the numerical computation of all quantities, we now choose $N = 64$ Fourier modes in space and a step size of $\tau = 2^{-13} T$ so that the homogenization error is not polluted by additional discretization errors. Figure \ref{error_plot} shows the homogenization errors over $\eta$ in double-logarithmic scaling for both the electric case (left) and the magnetic case (right) for which the observations are very similar. As highlighted by the reference lines, we observe the expected linear convergence in $\eta$ for the homogenized solution $\bfu_0$ and the expected quadratic convergence in $\eta$ for the first order corrections $\bfu_0 + \eta \bfu_1$. In addition, the second order corrections $\bfu_0 + \eta \bfu_1 + \eta^2 \bfu_2$ show the expected cubic convergence in $\eta$. Finally, we computed the second order macroscopic quantities $\bfu^{(2)}$ from sections \ref{sec:pure_marco_electric} and \ref{sec:pure_marco_magnetic}, which can be solved by the higher order wave equations \eqref{eq:pure_macro2_electric} and \eqref{eq:pure_macro2_magnetic}. Corrected with the microscopic quantities $\tbfu_1$ and $\tbfu_2$, the associated errors again show a cubic convergence in $\eta$, showing that $\bfu^{(2)}$ coincides with $\obfu^{(2)}$ up to higher order errors, but for which one has to solve only a single wave-type equation. \\
We also computed the errors of the macroscopic quantities, namely $\| \bfu_\eta - \bfu^{(1)} \|_{L^2(0,T;L^2)}$ and $\| \bfu_\eta - \bfu^{(2)} \|_{L^2(0,T;L^2)}$, which we do not show here. As expected, these errors do not exhibit higher-order decay but decrease only at first order in $\eta$. This indicates that the microscopic components play a crucial role in achieving higher-order approximation accuracy and cannot be neglected. However, although the order in $\eta$ is not improved by including the higher-order macroscopic components, we observed that the overall error magnitude is reduced. Hence, accounting for the nonlocal higher-order but purely macroscopic effects yields a more accurate approximation than the simple homogenized solution and can therefore be advantageous. \\
In summary, our numerical results verify our theoretical findings from Conclusion \ref{con:electric_hom}-\ref{con:magnetic_cor2} and the derived equations for the effective solution and its first and second order corrections lead to to an effective description of wave propagation through time-varying metamaterials.

\section{Conclusion}

In this work, we have developed a systematic higher-order homogenization framework for electromagnetic wave propagation in temporally modulated metamaterials with spatially homogeneous, time-dependent permittivity. Based on a formal two-scale asymptotic expansion, we derived the effective equations together with corrections up to second order, thereby capturing both local and nonlocal material laws. Our analysis reveals structural differences between the electric and magnetic field formulations and demonstrates how the position of the time-dependent coefficients determines the form of the resulting homogenized equations. A key contribution of this study is the careful derivation of consistent initial data and temporal interface conditions, which are intrinsic to the time-dependent setting and cannot be obtained through a straightforward space–time analogy. Our numerically verified formal results unify and extend previous approaches, providing a consistent mathematical framework for the effective modeling of time-varying metamaterials and offering a methodology that can be transferred and further developed for related wave systems, including acoustics, elasticity, and anisotropic metamaterials.

\subsubsection*{Acknowledgement}
This work is funded by the Deutsche Forschungsgemeinschaft (DFG, German Research Foundation) – Project-ID 258734477 – SFB 1173 and under Germany's Excellence Strategy – EXC-2047/1 – 390685813.

The authors would like to thank Puneet Garg, Michael Plum, Vishnu Raveendran and Carsten Rockstuhl for valuable discussions on the topic of this work.

\end{document}